\documentstyle[amssymb,amscd,12pt]{amsart}

\newtheorem{defn}{Definition}
\newtheorem{thm}[defn]{Theorem}

\newtheorem{cor}[defn]{Corollary}
\newtheorem{lem}[defn]{Lemma}
\newtheorem{prop}[defn]{Proposition}

\theoremstyle{remark}
\newtheorem{rem}{Remark}
\theoremstyle{remark}
\newtheorem{exam}{Example}

\numberwithin{equation}{section}
\numberwithin{defn}{section}

\begin{document}

%%%%%% 1.- OPERADORES
\newcommand\spk{{\operatorname{Spec}}(k)}
\renewcommand\sp{\operatorname{Spec}}
\renewcommand\sf{\operatorname{Spf}}
\newcommand\proj{\operatorname{Proj}}
\newcommand\aut{\operatorname{Aut}}
\newcommand\grv{{\operatorname{Gr}}(V)}
\newcommand\gr{\operatorname{Gr}}
\newcommand\glv{{\operatorname{Gl}}(V)}
\newcommand\glve{{\widetilde{\operatorname{Gl}}(V)}}
\newcommand\gl{\operatorname{Gl}}
\renewcommand\hom{\operatorname{Hom}}
\renewcommand\det{\operatorname{Det}}
\newcommand\detd{\operatorname{Det}^\ast}
\newcommand\im{\operatorname{Im}}
\newcommand\limi{\varinjlim}
\newcommand\limil[1]{\underset{#1}\varinjlim\,}
\newcommand\limp{\varprojlim}
\newcommand\limpl[1]{\underset{#1}\varprojlim\,}

%%%%%% 2.- CALIGRAFIA
\renewcommand\o{{\cal O}}
\renewcommand\L{{\cal L}}
\renewcommand\c{{\cal C}^\bullet}
\renewcommand\P{{\mathbb P}}
\newcommand\A{{\mathbb A}}
\newcommand\Z{{\mathbb Z}}

%%%%%% 3.- ABREVIATURAS
\newcommand\w{\widehat}
\renewcommand\tilde{\widetilde}

%%%%%% 4.- VARIOS
\newcommand\iso{@>{\sim}>>}
\renewcommand\lim{\underset{\underset{ A\sim V^+}A}\varprojlim}
\newcommand\M{{\cal M}}
\newcommand\fu{\underline}
\newcommand\kz{\fu{k((z))}^*}
\font\pun=cmsy10 at 5pt
\newcommand\punt {{\raise.4ex\hbox{\pun ^^O}}}
\newcommand\beq{
      \setcounter{equation}{\value{defn}}\addtocounter{defn}1
      \begin{equation}}

%%%%%%%%%%%%%%%%%%%%%%%%%%%%%%%%%%%%%%%%%%%%%%%%%%%%%%%%%%%%%
\title [The algebraic formalism of soliton equations] {The
algebraic formalism of soliton \\ equations over arbitrary base
fields}

\author[A. \'Alvarez \and J. M. Mu\~noz \and F. J. Plaza ] {A.
\'Alvarez V\'azquez \and \\ J. M. Mu\~noz Porras \and \\
F. J. Plaza Mart\'{\i}n \\ \medskip\tiny  Departamento de
Matem\'atica Pura y Aplicada \\ Universidad de Salamanca}

\address{Departamento de Matem\'atica Pura y Aplicada \\
Universidad de Salamanca \\ Plaza de la Merced 1-4 \\
 37008 Salamanca. Spain.}

\thanks{This work is partially supported by the CICYT research
contract n. PB91-0188 \\ \indent {\it Preprint:
}alg-geom/{\bf 9606009}}

\date{June 1996. (Update: November 1996, minor changes in \S5)}

\maketitle

%%%%%%%%%%%%%%%%%%%%%%%%%%%%%%%%%%%%%%%%%%%%%%%%%%%%%%%%%%%%%

\setcounter{tocdepth}1
\tableofcontents

%%%%%%%%%%%%%%%%%%%%%%%%%%%%%%%%%%%%%%%%%%%%%%%%%%%%%%%%%%%%%
\section{Introduction}

The aim of this paper is to offer an algebraic construction of
infinite-dimensional Grassmannians and determinant bundles.  As an
application we construct the $\tau$-function and formal
Baker-Akhiezer functions over arbitrary fields, by proving the
existence of a ``formal geometry'' of local curves analogous to the
geometry of global algebraic curves.

Recently G.~Anderson ([{\bf A}]) has constructed the
infinite-dimensional Grassmannians and
$\tau$-functions over $p$-adic fields; his constructions are
basically the same as in the Segal-Wilson paper ([{\bf SW}]) but he
replaces the use of the theory of determinants of Fredholm
operators over a Hilbert space by the theory of $p$-adic infinite
determinants (Serre [{\bf S}]).

Our point of view is completely different and the formalism used is
valid for arbitrary base fields; for example, for global number
fields or fields of positive characteristic.  We begin by defining
the functor of points, $\fu{\gr}(V,V^+)$, of the Grassmannian of a
$k$-vector space $V$ (with a fixed $k$-vector subspace
$V^+\subseteq V$) in such a way that the points
$\fu{\gr}(V,V^+)(\spk)$ are precisely the points of the
Grassmannian defined by Segal-Wilson or Sato-Sato ([{\bf SW}],
[{\bf SS}]) although the points over an arbitrary
$k$-scheme $S$ have been not previously considered by other
authors.  This definition of the functor $\fu{\gr}(V,V^+)$, which
is a sheaf in the category of $k$-schemes, allows us to prove that
it is representable by a separated $k$-scheme $\gr(V,V^+)$. The
universal property of the $k$-scheme $\gr(V,V^+)$ implies, as in
finite-dimensional Grassmannians, the existence of a universal
submodule,
$\L_V$, of $\pi^*V$ ($\pi:\gr(V,V^+)\to\spk$ being the natural
projection). These constructions allow us to use the theory of
determinants of Knudsen and Mumford ([{\bf KM}]) to construct the
determinant bundle over $\gr(V,V^+)$. This is one of the main
results of the paper because it implies that we can define
``infinite determinants'' in a completely algebraic way. From this
definition of the determinant bundle, we show in
$\S 3$ that global sections of the dual determinant bundle can be
computed in a very natural form. The construction of
$\tau$-functions and Baker functions is based on the algebraic
version, given in $\S4$, of the group $\Gamma$ of continuous maps
$S^1\to{\Bbb C}^*$ defined by Segal-Wilson ([{\bf SW}]) which acts
as a group of automorphisms of the Grassmannians. We replace the
group
$\Gamma$ by the representant of the following functor over the
category of
$k$-schemes
$$S\rightsquigarrow H^0(S,\o_S)((z))^*=
H^0(S,\o_S)[[z]][z^{-1}]^*$$ This is one of the points where our
view differs essentially from other known expositions ([{\bf A}],
[{\bf AD}], [{\bf SW}], [{\bf SS}]). Usually, the elements of
$\Gamma$ are described as developments, of the type
$f=\sum_{-\infty}^{+\infty}\lambda_i\,z^i$ ($\lambda_k\in{\Bbb
C}$), but in the present formalism the elements of $\Gamma$ with
values in a
$k$-algebra $A$ are developments
$f=\sum_{i\geq -N}^{}\lambda_i\,z^i\in A((z))$ such that
$\lambda_{-1},\dots,\lambda_{-N}$ are nilpotent elements of $A$.

In future papers we shall apply the formalism offered here to
arithmetic problems (Drinfeld moduli schemes and reprocity laws)
and shall give an algebraic formalism of the theory of KP-equations
related to the characterization of Jacobians and Prym varieties. We
also hope that this formalism might clarify the algebro-geometric
aspects of conformal field theories over base fields different from
$\Bbb R$ or $\Bbb C$ in the spirit of the paper of E.~Witten ([{\bf
W}]).

%%%%%%%%%%%%%%%%%%%%%%%%%%%%%%%%%%%%%%%%%%%%%%%%%%%%%%%%%%%%%
\section{Infinite Grassmannians\label{grass-section}}

Let $V$ be a vector space over a field $k$.

\begin{defn}{(Tate [{\bf T}])}  Two vector spaces $A$ and $B$ of
$V$ are commensurable if ${A+B}/{A\cap B}$ is a vector space over
$k$ of finite dimension. We shall use the symbol
$A\sim B$ to denote commensurable vector subspaces.
\end{defn}

Let us observe that commensurability is an equivalence relation
between vector subspaces. The addition and intersection of two
vector subspaces commensurable with a vector subspace $A$ is also
commensurable with $A$.

Let us fix a vector subspace $V^+\subseteq V$. The equivalence class
of vector subspaces commensurable with $V^+$ allows one to define
on $V$ a topology, which will be called
$V^+$-topology: a basis of neighbourhoods of $0$ in this topology
is the set of vector subspaces of $V$ commensurable with $V^+$.

$V$ is a Haussdorff topological space with respect to the
$V^+$-topology.

\begin{defn}
 The completion of $V$ with respect to the
$V^+$-topology is defined by:
$$ \w V= \underset{A\sim V^+}{\varprojlim} ( V/A)$$
\end{defn}

Analogously, given a vector subspace $B\subseteq V$ we define the
completions of $B$ and $ V/B$ with respect to $B\cap V^+$ and
${B+V^+}/{B}$,  respectively.

The homomorphism of completion $V@>i>> \w V$ is injective and $V$
is said to be complete if   $V@>i>> \w V$ is an isomorphism.

\begin{exam}
\begin{itemize}
\item $(V, V^+=0)$; $V$ is complete.
\item $V=k((t)),$  $V^+=k[[t]]$; $V$ is complete.
\item Let $(X, \o_X)$ be a smooth, proper and irreducible curve
over the field $k$, and let $V$ be the ring of adeles of the curve
and $V^+=\underset p{\prod} \w {\o_p}$ ( $\o_p$ being the $\frak
m_p$-adic completion of the local ring of $X$ in the point $p$);
$V$ is complete with respect the
$V^+$-topology.
\end{itemize}
\end{exam}

\begin{prop}  The following conditions are equivalent:
\begin{enumerate}
\item $V$ is complete.
\item $V^+$ is complete.
\item Each vector subspace commensurable with $V^+$ is complete.
\end{enumerate}
\end{prop}

\begin{pf}  This follows easily from the following commutative
diagram for every $A\sim V^+$:
$$ \CD 0 @>>> \widehat{A} @>>>  \widehat{V} @>>> \widehat{ V/A}
@>>> 0 \\ @. @A{i_A}AA @A{i_V}AA @A{\simeq}AA @. \\ 0  @>>>  A
@>>> V @>>>   V/A @>>> 0
\endCD$$
\end{pf}

\begin{defn}  Given a $k$-scheme $S$ and a vector subspace
$B\subseteq V$, we define:
\begin{enumerate}
\item $\w V_S=
 \lim ( V / A\underset k{\otimes}{\cal O}_S).$
\item $\widehat B_S=\lim ( {B} /{  A\cap B}) \underset
k{\otimes}{\cal O}_S$.
\item$\widehat{( V/B)}_S=\lim (( {V}/{A+B})\underset k{\otimes}
{\cal O}_S)$.
\end{enumerate}
\end{defn}

\begin{prop}
$\w{ V_S}$ is a sheaf of $\o_S$-modules and given $B\sim V^+$, we
have:
$$\w{(V/B)_S}={\w{ V_S}}/{\w{ B_S}}=(V/B)\underset k{\otimes} \o_S$$
\end{prop}

\begin{pf} This is an easy exercise of linear algebra.
\end{pf}

Let $V$ be a $k$-vector space and $V^+$ a vector subspace
determining a class of commensurable vector subspaces.

\begin{defn}
 A discrete vector subspace of $V$ is a vector subspace,
$L\subseteq V$, such that $L\cap V^+$ and
${V}/{L+V^+}$ are $k$-vector spaces of finite dimension.
\end{defn}

We aim to define a Grassmannian scheme $\gr(V,V^+)$, defining its
functor of points $\fu{\gr}(V,V^+)$ and proving that it is
representable in the category of $k$-schemes.

If $V$ is complete, the rational points  of our Grassmannian  will
be precisely the discrete vector spaces of $V$; that is,
$\fu{\gr}(V,V^+)(\spk)$ as a set coincides with the usual
 infinite Grassmannian defined by Pressley and Segal [{\bf PS}] or
M. and Y. Sato [{\bf SS}].

\begin{defn}
 Given a $k$-scheme $S$, a discrete submodule of $\w{V_S}$ is a
sheaf of quasi-coherent
$\o_S$-submodules $\L \subset \w{V_S}$ such that $\L_{S'}\subset
\w{V_{S'}}$ for every morphism $S'\to S$ and for each
$s\in S$, $\L\underset{\o_S}\otimes k(s)\subset
\w{V_S}\underset{\o_S}\otimes k(s)$ and  there exists an open
neighbourhood
$U_s$ of
$s$ and a commensurable $k$-vector subspace $B\sim V^+$ such that:
 $\L_{U_s}\cap \w{B_{U_s}}$ is free of finite type and
${\w{V_{U_s}}}/{\L_{U_s}+ \w{B_{U_s}}}=0.$
\end{defn}

\begin{prop} \label{prop-loc-free} With the notations of the above
definition, given another commensurable
$k$-vector space, $B'\sim V^+$, such that $B\subseteq B'$,
$\L_{U_s}\cap \w{{B'}_{U_s}}$ is locally free of finite type.
\end{prop}

\begin{pf}  This follows easily from the commutative diagram
$$
\CD 0 @>>> \w {B_{U_s}} @>>> \w  {{B'}_{U_s}} @>>>  (( B'/{B})
\underset k{\otimes} {\cal O}_{U_s})= {\w {B'_{U_s}}}/{\w
{B_{U_s}}} @>>> 0
\\ @. @VVV @VVV @VVV @. \\ 0 @>>>   {\w {V_{U_s}}}/{\L _{U_s}} @>>>
 {\w {V_{U_s}}}/{\L _{U_s}} @>>> 0 @>>> 0
\endCD$$   using the snake Lemma.
\end{pf}

\begin{defn}  Given a $k$-vector space $V$ and
$V^+\subseteq V$, the Grassmannian functor, $\fu{\gr}(V,V^+)$, is
the contravariant functor over the category of
$k$-schemes defined by
$$\fu{\gr}(V,V^+)(S)=\left\{
\begin{gathered} \text{discrete sub-$\o_S$-modules of $\w{V_S}$} \\
\text{with respect the $V^+$-topology} \end{gathered}
\right\}$$
\end{defn}

\begin{rem} Note that if $V$ is a finite dimensional
$k$-vector space and
$V^+=(0)$, then $\fu{\gr}(V,(0))$ is the usual Grassmannian functor
defined by Grothendieck [{\bf EGA}].
\end{rem}

\begin{defn}  Given a commensurable vector subspace
$A\sim V^+$, the functor $\fu{F_A}$ over the category of $k$-schemes
is defined by:
$$\fu{F_A}(S)=\{\text{ sub-$\o_S$-modules $\L\subset \w {V_S}$
such that $\L\oplus \w {A_S}=\w {V_S}$ }\}$$  (That is: $ \L \cap
\w {A_S}=(0)$ and  $ \L + \w {A_S}= \w {V_S}$).
\end{defn}

\begin{lem}\label{Frepre}  For every commensurable subspace $B\sim
V^+$, the contravariant functor $\fu{F_B}$ is representable by an
affine and integral $k$-scheme $F_B$.
\end{lem}

\begin{pf}  Let $L_0$ be a discrete $k$-subspace of $V$ such that
$L_0\oplus B = V$;  we then have:
$$ F_B(S)=\fu {Hom}_{\o_S}((\L_0)_S,\w B_S) =\lim
\fu{Hom}_{\o_S}((\L_0)_S,{ B}/{B\cap A}\underset
k{\otimes}\o_S))$$

If we denote by $F_{ B/{{B\cap A}}}(S)$ the set
 $\fu{ Hom}_{\o_S}((\L_0)_S, B/{B\cap A}{\underset k{\otimes}{\cal
O}_S)})$, it is obvious that the functor
$F_{ B/{{B\cap A}}}(S)$ is representable by an affine and integral
$k$-scheme since  $ B/{{B\cap A}}$ is a finite dimensional
$k$-vector space. But $\fu{F_B}$ is now a projective limit of
functors representable by affine schemes, so we conclude that
$\fu{F_B}$ is representable by an affine $k$-scheme.
\end{pf}

\begin{lem}\label{anterior} Let $\L$ be an element in
$\fu{\gr}_{V^+}(V)(S)$ and  $A$ and $B$ are two $k$-subspaces of
$V$ commensurable with $V^+$. It holds that:
\item{a)} if ${\w V_{S}} /{\L+ \w A_{S}} =0$, then
$\L\cap \w A_{S}$ is a finite type  locally free of $\o_S$-module.
\item{b)} ${\w V_{S}} /{\L+ \w B_{S}}$ is an
$\o_S$-module locally of finite presentation.
\end{lem}

\begin{pf}
\item{a)} By proposition {\ref{prop-loc-free}}, for each point
$s\in S$ there exists an open neighbourhood $U_s$ and a
commensurable $k$ subspace
$A'\sim V^+$ such that: $A\subseteq A'$,
$ { \w V_{U_s}}/{\L_{U_s}+
\w A'_{U_s}} =0$ and $\L_{U_s}\cap \w A'_{U_s}$ is free of finite
type.
  From the exact sequence:
$$0\to \L\cap \w A_{S} \to \L
\to ({\w V_S}/{\w A_{S}})= ( V/A)_S\to 0$$
 one deduces that $\L\cap \w A_{S}$ is quasicoherent  and
$$0\to (\L\cap \w A_{S})_{U_s}  \to \L_{U_s} \to ({\w V_{U_s}}/{\w
A_{U_s}})=( V/A)_{U_s} \to 0$$ Let us consider the commutative
diagram:
$$\CD 0 @>>> \L\oplus   \w {A_S} @>>>
\L\oplus  \w {{A'}_S } @>>>   ( {A'}/ A)_S  @>>> 0
\\ @. @VVV @VVV @VVV @. \\ 0 @>>> \w{V_S} @>>> \w {V_S}  @>>> 0
@>>> 0
\endCD$$  By using the snake lemma we have an exact sequence:
$$0\to (\L_{S}\cap \w A_{S})
\to \L_{S} \cap \w A'_{S} \to ( A'/A)_{S} \to  {\w{V}_{S}}
/{{\L_{S}}+ \w{ A }_{S}}\to  {\w{ V}_{S}}/{{\L_{S}}+ \w{
A'}_{S}}\to 0$$ In our conditions for $A$ and $A'$ we have:
$$0\to (\L_{U_s}\cap\w A_{U_s})
\to \L_{U_s} \cap \w A'_{U_s} \to ( A'/A)_{U_s} \to 0$$ Then,
$(\L\cap \w A_{S})_{U_s}=
\L_{U_s}\cap\w A_{U_s}$ is the kernel of a surjective homomorphism
between two free $\o_{U_s} $-modules of finite type, and we
conclude the proof.

\item{b)} For a given $s\in S$, let us take $B\subseteq B'$ such
that
$B'\sim V^+$ is a commensurable subspace,
 $\L_{U_s} \cap \w{B'}_{U_s}$ is free of finite type and
${\w{ V}_{U_s}} /{{\L_{U_s}}+ \w{ B'}_{U_s}}=0 $. We then have the
exact sequence:
$$  \L_{U_s} \cap \w{B'}_{U_s}\to
 B'/B \underset k\otimes \o_{U_s} \to
 {\w{ V}_{U_s}} /{{\L_{U_s}}+ \w{ B }_{U_s}} \to 0 $$ and so we
conclude.
\end{pf}

\begin{thm} The functor $\fu{\gr}(V,V^+)$ is representable by a
$k$-scheme
$\gr(V,V^+)$.
\end{thm}

\begin{pf}  The proof is modeled on the Grothendieck construction
of finite Grassmannians [{\bf EGA}]; that is:

It is sufficient to prove that $\{\fu{F_A}, A\sim V^+\}$ is a
covering of $\fu{\gr}(V,V^+)$ by open subfunctors:

\item {1)}{\sl For every $A\sim V^+$, the morphism of functors
$\fu{F_A}\to \fu{\gr}(V,V^+)$ is representable by an open
immersion}:
\newline That is, given morphism of functors
$X^{\punt}\to  \fu{\gr}(V,V^+)$ (where $X$ is  a $k$-scheme), the
functor
$$X^\punt \underset{\fu{\gr}(V,V^+)}\times \fu{F_A}
\hookrightarrow X^\punt $$ is represented by an open subscheme of
$X$. This is equivalent to proving that given $\L \in
\fu{\gr}(V,V^+)(X)$, the set :
$${\cal U}(A,\L)=\{x\in X\quad\text{ such that }\quad\L_x=
\L\underset{\o_X}\otimes k(x)\in F_A(\sp(k(x)))\}$$  is open in $X$.
\newline  If $\L_x  \in F_A(\sp(k(x))$, then:
$$  {\w V_{k(x)}}/{{\L_x} + \w A_{k(x)}} =0$$  but ${\w V_{X}}/{\L
+ \w {A_X}}$ is a $\o_X$-module of finite presentation and,
applying the lemma of Nakayama, there exists an open  neighbourhood
$U_x$ of
$x$, such that
$${\w V_{U_x}}/{\L_{U_x}+\w A_{U_x}}=0$$
\newline By lemma {\ref{anterior}}, $\L_{U_x}\cap \w {A_{U_x}}$  is
a
$\o_{U_x}$-module coherent. However, bearing in mind that   $\L_x
\cap \w {A_{k(x)} }=0$, there exists another open neighbourhood of
$x$, $U'_x\subseteq U_x$, such that  $\L_{U'_x}\cap \w
{A_{U'_x}}=0$ and therefore
$\L_{U'_x}\in F_A({U'}_x)$.

\item{2)}  {\sl For every $k$-scheme $X$ and every morphism of
functors
$$X^\punt\to\fu{\gr}(V,V^+)$$
 the open subschemes $\{ {\cal U}(A,\L),\, A\sim
V^+\}$ defined above are a covering of $X$.}
\newline That is, given $\L \in \fu{\gr}(V,V^+)(X)$ and a point
$x\in X$, there exists an open neighbourhood, $U_x$, of $x$ and a
commensurable subspace $A\sim V^+$ such that:
$$\L_{{U }_x}\in F_A({U }_x)$$     Let  $A$ be a commensurable
subspace such that:
$$ \L_x \cap \w {A_{k(x)}} =0$$  since  $ {\w V_{k(x)}}/{\L_x+ \w{
A_{k(x)}}}$ is a
$k(x)$-vector space of finite dimension, we can choose a basis
$(\overline {e_1 \otimes 1},
\dots,\overline {e_k\otimes 1})$  of $ {\w {V_{k(x)}}}/{  \L _x+ \w{
A_{k(x)}}}$ where
$e_i\in V$. Defining
$$B=A+\langle e_1,\dots,e_k\rangle$$ obviously $B\sim V^+$. One can
easily prove that there exists an open subset ${U'}_x \subseteq U$
such that
$\L_{{U'}_x} \in F_B({U'}_x)$ and  this completes the proof of the
theorem.
\end{pf}

\begin{lem}\label{anterior2} Let $A$, $B$ be two $k$-vector spaces
of $V$ commensurable with $V^+$. A necessary and sufficient
condition for the existence of $\L\in\fu{\gr}(V,V^+)(S)$ such that
 $\L\oplus\w{A_S}=\L\oplus\w{B_S}=\w{V_S}$, is that there should
exist an isomorphism of $k$-vector spaces
$$\tau: B/{A\cap B}\iso A /{A\cap B}$$
\end{lem}

\begin{pf} Let us consider the decomposition:
$$ \w {V_S}= \w {(A\cap B)_S}\oplus
\left( B /{A\cap B}\right)_S\oplus
\left(  { A /{A\cap B}}\right)_S\oplus \left( V/ {A+B}\right)_S$$
If the isomorphism $\tau$ exists, we take:
$$\L=\left\{(a,b,\tau(b)), \quad a \in \left( V/ {A+B}\right)_S
,\quad b\in \left( B/ {B\cap A}\right)_S \right\}$$  Conversely,
assume that  $\L\oplus\w{A_S}=\L\oplus\w {B_S}=\w{V_S}$. We then
have:
 $$\L \oplus \left( B /{B\cap A}\right)_S\oplus
\left(\w{B\cap A}\right)_S\simeq
\L \oplus  \w {B  _S}\simeq
\L\oplus  \w{ A _S}\simeq
\L \oplus \left( A/ {A\cap B}\right)_S \oplus
\left(\w {B\cap A}\right)_S$$
 from which we deduce that:
$$\left( B /{A\cap B}\right)_S\simeq
\left( A/{A\cap B}\right)_S$$
\end{pf}

\begin{thm}
$\gr(V,V^+)$ is a separated scheme.
\end{thm}

\begin{pf}  Let $F_B=\sp(A_B) =({\cal U}_B) $ be the affine open
subschemes of the Grassmannian constructed in lemma {\ref{Frepre}}.
It suffices to prove that given two commensurable subspaces $B'$
and $B$ such that $F_B\cap F_{B'}\neq\emptyset$ then $F_B\cap
F_{B'}$ is affine.

By lemma {\ref{anterior2}}
$$ F_B\cap F_{B'} \neq \emptyset$$ implies the existence of
$\L_0\in{F_B}(\spk)\cap F_{B'}(\spk)$ and bearing in mind that
$$ \fu{F_B} \cap\fu{F_{B'}} =
\fu{F_B} \underset{\fu{F_{B+B'}}}\times\fu{F_{B'}} $$ we conclude
the proof.
\end{pf}

\begin{defn}  The discrete submodule corresponding to the identity
$$Id\in \fu{\gr} (V,V^+)\left(\gr(V,V^+)\right)$$ will be called
the universal module and will be denoted by $$\L_V
\subset \w V_{\gr(V,V^+) }$$
\end{defn}

\begin{rem} In this section we have constructed infinite-dimensional
Grassmannian schemes in an abstract way. Since we select particular
vector spaces $(V,V^+)$ we obtain different classes of
Grassmannians. Two examples are relevant:
\begin{enumerate}
\item $V=k((t))$, $V^+=k[[t]]$. In this case, $\gr(k((t)),k[[t]])$
is the algebraic version of the Grassmannian constructed by
Pressley, Segal, and M. and Y. Sato ([{\bf PS}], [{\bf SS}]) and
this Grassmannian is particularly suitable for studying problems
related to the moduli of curves (over arbitrary fields) and
KP-equations.
\item Let $(X,\o_X)$ be a smooth, proper and irreducible curve over
the field $k$ and let $V$ be the adeles ring  over the curve and
$V^+=\underset p\prod
\w{\o_p}$ (Example 1.3.3). In this case  $\gr(V,V^+)$ is an adelic
Grassmannian which will be useful for studying arithmetic problems
over the curve $X$ or problems related to the classification of
vector bundles over a curve (non abelian theta functions...).
\newline Instead of adelic Grassmanians, we could define
Grassmanians associated with a fixed divisor on $X$ in an analogous
way.
\newline These adelic Grassmanians  will be also of interest in the
study of conformal field theories over Riemann surfaces in the
sense of Witten ([{\bf W}]).
\end{enumerate}
\end{rem}

\begin{rem} From the universal properties satisfied by the
Grassmannian one easily deduces the well known fact that given a
geometric point $W\in \gr(V)(Spec(K))$ ($k\hookrightarrow K$ being
an extension of fields), the Zariski tangent space to $\gr(V)$ at
the point $W$ is the
$K$-vector space:
$$T_W \gr(V)=Hom(W, \w{V_K}/W)$$
\end{rem}

%%%%%%%%%%%%%%%%%%%%%%%%%%%%%%%%%%%%%%%%%%%%%%%%%%%%%%%%%%%%%
%%% INDICE Y FIBRADO DETERMINANTE

\section{Determinant Bundles \label{det-bundles}}

In this section we construct the determinant bundle over the
Grassmannian following the idea of Knudsen and Mumford
([{\bf{KM}}]). This allow us to define determinants algebraically
and over arbitrary fields (for example for $k=\Bbb Q$ or $k={\Bbb
F}_q$).

Let us set a pair of vector spaces, $V^+\subseteq V$. As in section
{\ref{grass-section}}, we will denote the Grassmannian
$\gr(V,V^+)$ simply by $\grv$.

\begin{defn} For each $A\sim V^+$ and each $L\in\grv(S)$ we define
a complex,
$\c_A(L)$, of $\o_S$-modules by:
$$\c_A(L) \equiv \dots\to 0\to L\oplus \hat A_S  @>{\delta}>> \hat
V_S\to 0\to\dots$$
$\delta$ being the addition homomorphism.
\end{defn}

\begin{thm}
$\c_A(L)$ is a perfect complex of $\o_S$-modules.
\end{thm}

\begin{pf} We have to prove that the complex of $\c_A(L)$ is locally
quasi-isomorphic to a bounded complex of free finitely-generated
modules.

Let us note that the homomorphism of complexes given by the diagram:
 $$\CD
\dots @>>> 0 @>>> L\oplus \hat A_S @>{\delta}>>\hat  V_S @>>> 0
@>>> \dots
\\@. @. @V{p_1}VV @VVV @. @. \\
\dots@>>> 0 @>>> L @>{\phantom\delta}>> ( V/A)_S^{\hat{}} @>>> 0
@>>>\dots
\endCD$$ ($p_1$ being the natural projection) is a
quasi-isomorphism. The problem is local on $S$, and hence for each
$s\in S$ we can assume the existence of an open neighbourhood, $U$,
and a commensurable subspace $B\sim V^+$ such that $A\subseteq B$
and:
$$ \hat V_U/{(L_U,\hat B_U)} = 0\qquad,\qquad  L_U\cap \hat B_U
\text{ is free and finitely-generated}$$ We then have the exact
sequence:
$$0 \to L_U\cap\hat  A_U \to L_U\cap\hat  B_U \to ( B/A)_U \to{\hat
V_U}/{(L_U+\hat A_U)} \to 0$$ from which we deduce that the
homomorphism of complexes given by the following diagram is a
quasi-isomorphism:
$$\CD
\dots @>>> 0 @>>> L_U\cap\hat B_U  @>>>( B/A)^{\hat{}}_U @>>> 0
@>>>\dots
\\ @. @. @VVV @VVV @. @.\\
\dots @>>> 0 @>>> L_U  @>>> ( V/A)^{\hat{}}_U @>>> 0 @>>>\dots
\endCD$$ That is, $\c_A(L)\vert_U$ is quasi-isomorphic to the
complex
$0 \to L_U\cap\hat B_U \to ( B/A)_U \to 0$, which is a complex of
free and finitely-generated modules.
\end{pf}

\begin{defn} The index of a point $L\in\grv(S)$ is the locally
constant function
$i_L\colon S\to\Z$ defined by:
$$i_L(s)=\text{  Euler-Poincar\'e characteristic of }
\c_{V^+}(L)\otimes k(s)$$
$k(s)$ being the residual field of the point $s\in S$. (For the
definition of the Euler-Poincar\'e characteristic of a perfect
complex see {\rm [{\bf{KM}}]}).
\end{defn}

\begin{rem} The following properties of the index are easy to
verify:
\begin{enumerate}
\item  Let $f:T\to S$ be a morphism of schemes and $L\in\grv(S)$;
then: $i_{f^*L}=f^*(i_L)$.
\item The function $i$ is constantly zero over the open subset
$F_{V^+}$.
\item If $B\sim V^+$, $V^+\subseteq B$ and ${\hat V_U}/{L_U+\hat
B_U}=0$ over an open subscheme $U\subseteq S$, then
$i_L(s)=\dim_{k(s)}(L_s\cap \hat B_s)-\dim_{k(s)}({B_s}/{V^+_s})$.
\item If $V$ is a finite-dimensional $k$-vector space and $V^+=V$
and $L\in\grv(S)$, then $i_L=\operatorname{rank}(L)$.
\item For any rational point $L\in\grv(\spk)$ one has:
$$i_L=\dim_k(L\cap V^+)-\dim_k({\hat V}/{L+\hat V^+})$$
\end{enumerate}
\end{rem}

\begin{thm} Let $\gr^n(V)$ be the subset over which the index takes
values equal to $n\in\Z$. $\gr^n(V)$ are open connected
subschemes of
$\grv$ and the decomposition of $\grv$ in connected components is:
$$\grv=\underset{n\in\Z}\coprod\gr^n(V)$$
\end{thm}

\begin{pf} This is obvious from the properties of the index.
\end{pf}

Given a point $L\in\grv(S)$ and $A\sim V^+$, we denote by
$\det\c_A(L)$ the determinant sheaf of the perfect complex
$\c_A(L)$ in the sense of [{\bf{KM}}].

\begin{thm} With the above notations the invertible sheaf over $S$,
$\det\c_A(L)$, does not depend on $A$ (up to isomorphisms).
\end{thm}

\begin{pf} Let $A$ and $A'$ be two commen\-surable subspaces. It
suffices to prove that:
$$\det\c_A(L)\iso\det\c_{A'}(L)$$  in the case $A\subseteq A'$. In
this case we have a diagram:
$$\CD
\dots@>>> 0@>>> \L\oplus\hat A_S @>{\delta}>>\hat V_S @>>>0@>>>\dots
\\ @. @. @VVV @V{Id}VV @. @. \\
\dots@>>>0@>>>\L\oplus\hat A'_S @>{\delta}>>\hat V_S @>>>0@>>>\dots
\endCD$$ and by the additivity of the functor $\det(-)$ we obtain:
$$\det\c_A(L)\otimes\det( {A'}/A)_S\iso \det\c_{A'}(L)$$ However $
{A'}/A$ is free and we conclude the proof.
\end{pf}

\begin{defn} The determinant bundle over $\grv$, $\det_V$, is the
invertible sheaf: $$\det\c_{V^+}(\L_V)$$
 $\L_V$ being the universal submodule over $\grv$.
\end{defn}

\begin{prop}{(Functoriality)} Let $L\in\grv(S)$ be  a point given
by a morphism $f_L:S\to\grv$. There exists a functorial isomorphism:
$$f^*_L\det_V\iso\det\c_A(L)$$ We shall denote this sheaf by
$\det_V(L)$.
\end{prop}

\begin{pf} The functor $\det(-)$ is stable under base changes.
\end{pf}

\begin{rem} Let $L\in\grv(\sp(K))$ be a rational point and let
$A\sim V^+$ such that $L\cap \hat A_K$ and ${\hat V_K}/{L+\hat
A_K}$ are
$K$-vector spaces  of finite dimension. In this case we have an
isomorphism:
$$\det_V(L)\simeq
\wedge^{max}(L\cap\hat A_K)\otimes\wedge^{max}( {\hat V_K}/{(L+\hat
A_K)})^\ast$$ That is, our determinant coincides, over the
geometric points, with the determinant bundles of Pressley, Segal,
Wilson and M. and Y. Sato ([{\bf PS}], [{\bf SW}], [{\bf SS}]).
\end{rem}

We shall now state with precision the connection between the
determinant bundle $\det_V$ and the determinant bundle over the
finite Grassmannianns.

Let $L,L'\in\grv(\spk)$ such that $L\subseteq L'$. In these
conditions,
${L'}/L$ is a $k$-vector space of finite dimension. The natural
projection $\pi:L'\to{L'}/L$ induces an injective morphism of
functors:
$$\fu\gr({L'}/L)\hookrightarrow\fu\gr(V)$$ defined by:
$$j(M)=\pi^{-1}(M)\qquad\text{for each }M\in\fu\gr({L'}/L)(S)$$ We
then have a morphism of schemes:
$$j:\gr({L'}/L)\hookrightarrow\grv$$ It is not difficult to prove
that $j$ is a closed immersion.

\begin{thm} With the above notations, there exists a natural
isomorphism:
$$j^*\det_V\iso\det_{{L'}/L}$$
$\det_{{L'}/L}$ being the determinant bundle over the finite
Grassmannian $\gr({L'}/L)$.
\end{thm}

\begin{pf}
Let $\L^f$ be the universal submodule over
$\gr({L'}/L)$. By definition $\det_{{L'}/L}=\det(\L^f\to{L'}/L)$,
which is isomorphic to $\det(\pi^{-1}\L^f\to L')$. By the
definition of $j$, one has $j^*\L_V\iso\pi^{-1}\L^f$ and hence:
$$j^*\det_V\simeq\det(\pi^{-1}\L^f\oplus\hat V^+\to\hat V)$$ and
from the exact sequence of complexes:
$$\CD
\pi^{-1}\L^f @>>> \pi^{-1}\L^f\oplus \hat V^+ @>>> \hat V^+
\\ @VVV @VVV @VVV \\ L' @>>> \hat V @>>> {\hat V}/{L'}
\endCD$$ we deduce that $j^*\det_V\simeq\det_{{L'}/L}$.
\end{pf}

\begin{cor} Let $i$ be the index function over $\grv$. For each
rational point
$M\in\gr({L'}/L)$ one has:
$$i(j(M))=i(L')+\dim_k({L'}/{M+L})$$
\end{cor}

\begin{pf} Obvious.
\end{pf}

%%%%%%%%%%%%%%%%%%%%%%%%%%%%%%%%%%%%%%%%%%%%%%%%%%%%%%%%%%%%%
%%% SECTIONS OF THE DETERMINANT BUNDLE

\subsection{Global sections of the determinant bundles and
Pl\"ucker morphisms \label{det-section}}

It is well known that the determinant bundle have no global
sections. We shall therefore explicitly construct global sections
of the dual of the determinant bundle over the connected component
$\gr^0(V)$ of index zero.

We use the following notations: $\wedge^\bullet E$ is the exterior
algebra of a $k$-vector space  $E$; $\wedge^r E$ its component of
degree $r$, and $\wedge E$ is the component of higher degree when
$E$ is finite-dimensional.

Given a perfect complex $\c$ over  $k$-scheme $X$, we shall write
$\det^*\c$ to denote the dual of the invertible sheaf $\det\c$.

To explain how global sections of the invertible sheaf $\det^*\c$
can be constructed, let us begin with a very simple example:

\noindent Let $f:E\to F$ be a homomorphism between
finite-dimensional $k$-vector spaces of equal dimension. This
homomorphism induces:
$$\wedge(f):\wedge E\to\wedge F$$ and $\wedge(f)\ne 0 \iff f$ is an
isomorphism.
$\wedge(f)$ can be expressed as a homomorphism:
$$\wedge(f):k\to \wedge F\otimes (\wedge E)^*$$ Thus, if we
consider $E@>{f}>> F$ as a perfect complex, $\c$, over $\spk$, we
have defined a {\bf canonical section}
$\wedge(f)\in H^0(\spk,\det^*\c)$.

Let us now consider a perfect complex
$\c\equiv(E@>{f}>>F)$ of sheaves of $\o_X$-modules over a
$k$-scheme $X$, with Euler-Poincar\'e characteristic
${\cal X}(\c)=0$. Let $U$ be an open subscheme of $X$ over which
$\c$ is quasi-isomorphic to a complex of finitely-generated free
modules. By the above argument, we construct a canonical section
$det(f\vert_U)\in H^0(U,\det^*\c)$ and for other open subset,
$V$, there is a canonical isomorphism $det(f\vert_U)\vert_{U\cap
V}\simeq det(f\vert_V)\vert_{U\cap V}$ and we therefore  have a
canonical section $\det(f)\in H^0(X,\det^*\c)$. If the complex
$\c$ is acyclic, one has an isomorphism:
$$\aligned
\o_X &\iso\det^\ast\c\\ 1 &\mapsto det(f)
\endaligned$$ (for details see [{\bf KM}]).

Let $0\to H^\bullet\to\c_1\to\c_2\to 0$ be an exact sequence of
perfect complexes. There exists a functorial isomorphism
$$\det^\ast\c_1\iso \det^\ast
H^\bullet\otimes_{\o_X}\det^\ast\c_2$$  If $H^\bullet$ is acyclic,
we obtain an isomorphism
$$ H^0(X,\det^\ast\c_2)\iso
 H^0(X,\det^\ast\c_1)$$ In the case $H^\bullet\equiv (E @>{Id}>>
E)$,
$\c_1\equiv (V@>{f}>>V)$, $\c_2\equiv (F @>{f'}>>F)$, and ${\cal
X}(\c_i)=0$ ($i=1,2$), we obtain the following commutative diagram:
$$\CD
\o_X @>{\sim}>>\det^\ast H^\bullet\otimes\o_X \simeq\o_X
\\@VVV @VVV \\
\det^\ast\c_1   @>{\sim}>>
\det^\ast H^\bullet\otimes\det^\ast\c_2\simeq\det^\ast\c_2
\endCD$$ from which we deduce that $det(f)=det(f')$. Moreover, if
$F$ is locally free of finite rank, this means that  computation of
$det(f)$ is reduced to computation of
$det(f')$, as mentioned above.

Let $V$ be a $k$-vector space and $V^+\subseteq V$ and $A\sim V^+$
a commensurable vector subspace. Let us consider the perfect
complex $\c_A\equiv(\L\oplus \hat A@>{\delta_A}>>\hat V)$ over
$\grv$ defined in {\ref{det-bundles}} ($\L$ being the universal
discrete submodule over $\grv$).

\begin{lem}
$F_A\subseteq\gr^0(V)$ if and only if
$$\dim_k( A/{A\cap V^+})-\dim_k({V^+}/{A\cap V^+})=0$$
\end{lem}

\begin{pf} Obvious.
\end{pf}

\begin{cor} The open subschemes $F_A$ with $\dim_k( A/{A\cap
V^+})-\dim_k({V^+}/{A\cap V^+})=0$ are a covering of
$\gr^0(V)$. Given $A,B\sim V^+$ under the assumption
$F_A,F_B\subseteq\gr^0(V)$ one has $\dim_k( A/{A\cap B})-\dim_k(
B/{A\cap B})=0$.
\end{cor}

\begin{pf} Obvious.
\end{pf}

Given $A\sim V^+$ with $F_A\subseteq \gr^0(V)$, let us note that
$\c_A\vert_{F_A}$ is an acyclic complex. One then has an
isomorphism:
$$\aligned
\o_X\vert_{F_A}&@>\sim>>
\det^\ast\c_A\vert_{F_A}\\  1 &\mapsto s_A=det(\delta_A)\vert_{F_A}
\endaligned$$ We shall prove that the section
$s_A\in H^0(F_A,\det^*\c_A)$ can be extended in a canonical way
to a global section of $\det^*\c_A$ over the Grassmannian
$\gr^0(V)$.

Let $B\sim V^+$ be such that $F_B\subseteq\gr^0(V)$ and let us
consider the complex:
$$\c_{AB}\equiv(\L\oplus \hat A@>{\delta_{AB}}>>
\L\oplus \hat B)$$ where $\delta_{AB}=\delta_B^{-1}\circ\delta_A$.
Obviously $\delta_{AB}\vert_{(0,\hat A\cap\hat B)}= Id_{\hat
A\cap\hat B}$ and
$\delta_{AB}\vert_{(\L,0)}= Id_\L$, we then have an exact sequence
of complexes:
$$\CD
 @. @. @.  A/{A\cap B}
\\ @. @. @. @V{\simeq}VV \\ 0@>>>\L\oplus(\hat A\cap\hat B)
@>>>\L\oplus\hat A @>>>
 {(\L\oplus\hat A)}/{\L\oplus(\hat A\cap\hat B)}@>>>0
\\ @. @V{Id}VV @V{\delta_{AB}}VV @V{\phi_{AB}}VV @. \\ 0
@>>>\L\oplus(\hat A\cap\hat B) @>>>\L\oplus\hat B @>>>
{(\L\oplus\hat B)}/{\L\oplus(\hat A\cap\hat B)}@>>> 0
\\ @. @. @. @V{\simeq}VV \\ @. @. @.  B/{A\cap B} \\
\endCD$$ and from the discussion at the beginning of this section
we have that $det(\phi_{AB})=det(\delta_{AB})\in
 H^0(F_B,\det^\ast\c_{AB})$ and $det(\delta_{AB})$ satisfies the
cocycle condition:
$$\aligned &det(\delta_{AA})=1 \\ &det(\delta_{AB})\cdot
det(\delta_{BC}) = det(\delta_{AC})
\qquad\text{over $F_B\cap F_C$ for any }C\sim V^+
\endaligned$$

Over $F_A\cap F_B$ we have canonical isomorphisms:
$$\aligned
\o_{F_A\cap F_B}&\iso \det^\ast\c_A\vert_{F_A\cap F_B} \\ 1&\mapsto
s_A \endaligned$$
$$\aligned
\o_{F_A\cap F_B}&\iso \det^\ast\c_B\vert_{F_A\cap F_B} \\ 1&\mapsto
s_B \endaligned$$
$$\aligned
\o_{F_A\cap F_B}&\iso \det^\ast\c_{AB}\vert_{F_A\cap F_B} \\
1&\mapsto det(\delta_{AB}) \endaligned$$ which are compatible,
therefore:
$$(s_B\cdot det(\delta_{AB}))\vert_{F_A\cap F_B}= s_A\vert_{F_A\cap
F_B}$$
$s_B\cdot det(\delta_{AB})$ being the image of $s_B\otimes
det(\delta_{AB})$ by the homomorphism:
$$ H^0(F_B,\det^\ast\c_B)\otimes
 H^0(F_B,\det^\ast\c_{AB})\to\
 H^0(F_B,\det^\ast\c_A)$$ defined by the isomorphism of sheaves:
$$\det\c_A\simeq \det\c_B\otimes \wedge( A/{A\cap B})\otimes
\wedge( B/{A\cap B})^\ast\simeq
\det\c_B\otimes\det\c_{AB}$$

\begin{defn} The global section
$\omega_A\in H^0(\gr^0(V),\det^*\c_A)$ defined by:
$$\{s_B\cdot det(\delta_{AB})\}_{B\sim V^+}$$
 will be called the canonical section of $\det^\ast\c_A$.
\end{defn}

This result allows us to compute many global sections of
$\detd_V=\det\c_{V^+}$ over $\gr^0(V)$:

\noindent Given $A\sim V^+$ such that $F_A\subseteq\gr^0(V)$ the
isomorphism $\det^*\c_A\iso\detd_V$ is not canonical, and in fact
we have a canonical isomorphism:
$$\det^*\c_A\iso\detd_V\otimes\bigwedge( A/{A\cap V^+})\otimes
\bigwedge({V^+}/{A\cap V^+})^*$$ Therefore to give an isomorphism
$\det^*\c_A\iso\detd_V$ depends on the choice of bases for the
vector spaces
$ A/{A\cap V^+}$ and ${V^+}/{A\cap V^+}$.

%%%%%%%%%%%% FINITE DIMENSIONAL CASE
\subsection{Computations for finite-dimensional
Grassmannians.}\label{comp-fin-dim}

Let $V$ be a d-dimensional $k$-vector space with a basis
$\{e_1,\dots,e_d\}$, let $\{e^\ast_1,\dots,e^\ast_d\}$ be its dual
basis, and $V^+=<e_{k+1},\dots,e_d>\subseteq V$. In
this case $\gr^0(V,V^+)$ is the Grassmannian of $V$
classifying $k$-dimensional vector subspaces of $V$. Given a
family of indexes $1\leq i_1<\dots< i_l\leq d$ ($1\le l\le d$), let
$A(i_1,\dots,i_l)$ be the vector subspace generated by
$\{e_{i_1},\dots,e_{i_l}\}$. One has that $F_A\subseteq\gr^0(V)$ is
equivalent to saying that $l=d-k$.

Let us set $A=A(i_1,\dots,i_{d-k})$. Now, the canonical section
$\omega_A\in H^0(\gr^0(V),\det^\ast\c_A)$ is the section whose
value at the point $L=<l_1,\dots,l_k>\in\gr^0(V)$ is given by:
$$\omega_A(L)=\pi_A(l_1)\wedge\dots\wedge\pi_A(l_k)
\otimes l^\ast_1\wedge\dots\wedge l^\ast_k
\in\wedge V/A\otimes\wedge L^\ast=(\det^\ast\c_A)_L$$
$\{l^\ast_1,\dots,l^\ast_k\}$ being the dual basis of
$\{l_1,\dots,l_k\}$ and
$\pi_A\colon L\to V/A$ the natural projection. Note that
$\{e_{j_1},\dots,e_{j_k}\}$ is a basis of $V/A$ where
$\{j_1,\dots,j_k\}=\{1,\dots,d\}-\{i_1,\dots,i_{d-k}\}$, and
that its dual basis is $\{e^\ast_{j_1},\dots,e^\ast_{j_k}\}$ in $(
V/A)^\ast\subset V^\ast$. We have now:
$$\omega_A(L)= (e^\ast_{j_1}\wedge\dots\wedge e^\ast_{j_k})
(l_1\wedge \dots\wedge l_k)\cdot  e_{j_1}\wedge\dots\wedge
e_{j_k}\otimes l^\ast_1\wedge \dots\wedge l^\ast_k$$

Observe that the $k$-vector space
$$\wedge( A/{A\cap V^+})^*\otimes\wedge( V^+/{A\cap V^+})$$
is generated by:
$$e_A=e^*_{m_1}\wedge\dots\wedge e^*_{m_r}\otimes
e_{n_1}\wedge\dots\wedge e_{n_r}$$
where
$$\aligned
&\{i_1,\dots,i_{d-k}\}-\{k+1,\dots,d\}=\{m_1,\dots,m_r\}\\
&\{k+1,\dots,d\}-\{i_1,\dots,i_{d-k}\}=\{n_1,\dots,n_r\}
\endaligned$$ And tensorializing by $e_A$ gives an isomorphism:
$$ H^0(\gr^0(V),\det^\ast\c_A) @>{\otimes e_A}>>
 H^0(\gr^0(V),\detd_V)$$
Let $\Omega_A$ be the image of the canonical section
$\omega_A$. The explicit expression of $\Omega_A$ is:
$$\Omega_A(L)=
(e^*_{j_1}\wedge\dots\wedge e^*_{j_k})(l_1\wedge \dots\wedge l_k)
\cdot
e_1\wedge\dots\wedge e_{d-k}\otimes l^*_1\wedge \dots\wedge l^*_k
\in(\det_V^\ast)_L$$

Let $\L\subseteq V_{\gr^0(V)}$ be the universal submodule. One has
a canonical epimorphism:
$$\wedge^k V^*_{\gr^0(V)}\to\wedge^k\L^*$$
 and bearing in mind the canonical isomorphism
$\detd_V\simeq \wedge( V/{V^+})\otimes\wedge^k\L^*$ one obtains a
canonical homomorphism:
$$\aligned
\wedge^k V^\ast= H^0(\gr^0(V),\wedge^k V^\ast)&\to
 H^0(\gr^0(V),\detd_V)\\
e^*_{j_1}\wedge\dots\wedge e^*_{j_k}&\longmapsto
\Omega_{A(i_1,\dots,i_{d-k})}
\endaligned$$
(where
$\{j_1,\dots,j_k\}\coprod\{i_1,\dots,i_{d-k}\}=\{1,\dots,d\}$).

It is well known that this homomorphism is in fact an isomorphism.

%%%%%%%%% INFINITE-DIMENSIONAL CASE
\subsection{Computations for infinite-dimensional Grassmannians.}

In {\ref{comp-fin-dim}} we discussed well known facts about the
determinants of finite-dimen\-sional Grassmannians but have stated
these results in an intrinsic language, which can easily be
generalized to the infinite-dimensional case.

Let $V$ be a $k$-vector space. We shall assume that there exists a
family of linearly independent vectors
$\{e_i,i\in\Z\}$ such that:
\begin{enumerate}
\item $<\{e_i\}, i\ge 0>$ is dense in $\hat V^+$ (with
respect to the $V^+$-topology),
\item $<\{e_i\}, i\in\Z>$ is dense in $\hat V$.
\end{enumerate}

\begin{rem} The above conditions are satisfied for example by
$V=k((t))$ and $V^+=k[[t]]$.
\end{rem}

\begin{defn} Let ${\cal S}$ be the set of sequences
$\{s_0,s_1,\dots\}$ of integer numbers satisfying the following
conditions:
\begin{enumerate}
\item the sequence is strictly increasing,
\item there exists $s\in\Z$ such that
$\{s,s+1,s+2,\dots\}\subseteq\{s_0,s_1,\dots\}$,
\item $\#(\{s_0,s_1,\dots\}-\{0,1,\dots\})=
\#(\{0,1,\dots\}-\{s_0,s_1,\dots\})$.
\end{enumerate}

\end{defn}

The sequences of ${\cal S}$ are usually called Maya's diagrams or
Ferrer's diagrams of virtual cardinal zero (this is condition 3).

For each $S\in{\cal S}$, let $A_S$ be the vector subspace of $V$
generated by $\{e_{s_i}, i\ge 0\}$. By the condition 3 one has:
$$\dim_k({A_S}/{A_S\cap V^+})=
\dim_k({V^+}/{A_S\cap V^+})$$ and hence: $A_S\sim V^+$ and
$F_{A_S}\subseteq\gr^0(V)$. Further, $\{F_{A_S},S\in{\cal S}\}$ is
a covering of
$\gr^0(V)$.

Let $\{e^*_i\}$ be a dual basis of $\{e_i\}$; that is, elements of
$V^*$ given by $e^*_i(e_j)=\delta_{ij}$.

For each finite set of increasing integers,
$J=\{j_1,\dots,j_r\}$, let us define
$e_J=e_{j_1}\wedge\dots\wedge e_{j_r}$ and
$e^*_J=e^*_{j_1}\wedge\dots\wedge e^*_{j_r}$.

Given $S\in{\cal S}$, choose $J,K\subseteq\Z$ such that
$\{e_j\}_{j\in J}$ is a basis of ${A_S}/{A_S\cap V^+}$ and
$\{e^*_k\}_{k\in K}$ of ${V^+}/{A_S\cap V^+}$. We have seen that
tensorializing by $e_J\otimes e^*_K$ defines an isomorphism:
$$ H^0(\gr^0,\det^*\c_{A_S})  @>{\,\otimes(e_J\otimes e^*_K)\,}>>
 H^0(\gr^0,\detd_V)$$

\begin{defn}\label{global-section} For each $S\in{\cal S}$,
$\Omega_S$ is the global section of $\detd_V$ defined by:
$$\Omega_S=\omega_{A_S}\otimes e_J\otimes e^*_K$$ We shall denote
by $\Omega_+$ the canonical section of
$\detd_V$.
\end{defn}

Let $\Omega({\cal S})$ be the $k$-vector subspace of
$ H^0(\gr^0,\detd_V)$ generated by the global sections
$\{\Omega_S,S\in{\cal S}\}$.

We define the Pl\"ucker morphism:
$$\aligned{\cal P}_V: \gr^0(V) &\to \check\P\Omega(S) \\ L &\mapsto
\{\Omega_S(L)\}\endaligned$$ as the morphism of schemes defined by
the homomorphism of sheaves:
$$\Omega(S)_{\grv}\to\detd_V\to 0$$
(by the universal property of $\check\P$).

\begin{rem}
Given $L,L'\in\grv(\spk)$ such that $L\subseteq L'$,
let  $j:\gr^0({L'}/L)\hookrightarrow\gr^0(V)$ be the natural closed
immersion. Since $j^*\det_V\simeq\det_{{L'}/L}$, one can easily see
that the composition:
$$ \gr^0({L'}/L)\overset j\hookrightarrow\gr^0(V)
\overset{\quad{\frak p}\quad}\to \check\P\Omega(S)$$ factors
through the Pl\"ucker immersion of the finite-dimensional
Grassmannian $\gr^0({L'}/L)$:
$${\frak p}_{{L'}/L}:\gr^0({L'}/L)\to
\proj S^\punt H^0(\gr^0({L'}/L),\det^*_{{L'}/L})$$
\end{rem}

\begin{thm}
The Pl\"ucker morphism is a closed immersion.
\end{thm}

\begin{pf}
Going on with the analogy with finite grassmannians, we will show
that this morphism is locally given as the graph of a suitable
morphism. Consider the morphism:
$$F_{A_S}\hookrightarrow\gr^0(V) @>{\cal P}>> \check\P\Omega(S)$$
From the universal property of $\check\P$, we deduce a epimorphism:
$$f_S:\Omega(S)\underset{k}\otimes B\to B$$
(where $\sp(B)=F_{A_S}$, and $\detd_V\vert_{F_{A_S}}$ is a line
bundle). Note that it has a section, since the image of $\Omega_S$ is
a everywhere non null function. That is, there exists a subspace
$W\subset \Omega(S)$, and an isomorphism of $k$-vector spaces:
$$<\Omega_S>\oplus W \iso \Omega(S)$$
such that $f_S$ is the projection onto the first factor. In other
words, ${\cal P}\vert_{F_{A_S}}$ is the graph of a morphism.
\end{pf}

\begin{rem}\label{section-ogr}
Note that considering the chain of finite-dimensional Grassmannians
$\gr^0({L_i}/{L_{-i}})$ ($L_i$ being the subspaces $<\{e_j\}_{j\leq
i}>$),  which are closed subschemes of $\gr^0(V)$, one easily
deduces that $ H^0(\gr^0(V),\o_{\gr^0(V)})=k$ from the fact that
the homomorphism:
$$ H^0(\gr^0(V),\o_{\gr^0(V)})\to\limp
 H^0(\gr^0({L_i}/{L_{-i}}),\o_{\gr^0({L_i}/{L_{-i}})})$$ is
injective.
\end{rem}

%%%%%%%%%%%%%%%%%%%%%%%%%%%%%%%%%%%%%%%%%%%%%%%%%%%%%%%%%%%%%
%%% AUTOMORPHISMS OF THE GRASSMANNIAN

\section{Automorphisms of the Grassmannian and the
 ``formal geometry'' of local curves}
\label{aut-grass}

Let $(V,V^+)$ be a pair of a $k$-vector space and a vector subspace
$V^+\subseteq V$ and let $\grv$ denote the corresponding
Grassmannian. We shall define the algebraic analogue of the
restricted linear group defined by Pressley, Segal and Wilson
([{\bf PS}], [{\bf SW}]). This group is too large to be
representable by a $k$-scheme and we therefore define it as a sheaf
of groups in the category of
$k$-schemes.

For each $k$-scheme $S$, let us denote by
$\aut_{\o_S}(\hat V_S)$ the group of automorphisms of the
$\o_S$-module $\hat V_S$.

\begin{defn}
\item{a)} A sub-$\o_S$-module ${\cal B}\subseteq \hat V_S$ is said
to be locally commensurable with $V^+$ if for each $s\in S$ there
exists an open neighbourhood $U_s$ of $s$ and a commensurable
vector subspace $B\sim V^+$ such that ${\cal B}\vert_{U_s}=\hat
B_{U_s}$.
\item{b)} An automorphism $g\in\aut_{\o_S}(\hat V_S)$ is called
bicontinuous with respect to the $V^+$-topo\-logy if
$g(\hat V^+_S)$ and $g^{-1}(\hat V^+_S)$ are
$\o_S$-modules of $\hat V_S$ locally commensurable with
$V^+$.
\item{c)} The linear group, $\glv$, of $(V,V^+)$ is the
contravariant functor over the category of $k$-schemes defined by:
$$S\rightsquigarrow \glv(S)=\{g\in\aut_{\o_S} (\hat V_S)\text{ such
that $g$ is bicontinuous }\}$$
\end{defn}

\begin{thm} There exists a natural action of $\glv$ over the functor
of points of the Grassmannian $\grv$:
$$\aligned
\glv\times & \fu\grv @>{\mu}>>\fu\grv  \\ (g, &
L)\phantom{xx}\longmapsto g(L)
\endaligned$$
\end{thm}

\begin{pf} Let $g\in\glv(S)$ and  $L\in\fu\grv(S)$. We have:
$${\hat V_S}/{g(L)+\hat V^+_S}\simeq  {\hat V_S}/{L+g^{-1}{\hat
V^+_S}}$$ and by definition of bicontinuous automorphisms, for each
$s\in S$ there exist an open neighbourhood $U_s$ and a
commensurable $A\sim V^+$ such that $g^{-1}\hat V^+\vert_{U_s}=\hat
A_{U_s}$. Then:
$${\hat V_{U_s}}/{g(L)_{U_s}+\hat V^+_{U_s}}\simeq  {\hat
V_{U_s}}/{L_{U_s}+\hat A^+_{U_s}}$$ from which we deduce that
$g(L)\in\fu\grv(S)$.
\end{pf}

\begin{thm} There exists a canonical central extension of functors
of groups over the category of $k$-schemes:
$$0\to{\Bbb G}_m\to \glve \to \glv \to 0$$ and a natural action
$\bar\mu$ of $\glve$ over the vector bundle
${\Bbb V}(\det_V)$ defined by the determinant bundle, such that the
following diagram is commutative:
$$\CD
\glve\times{\Bbb V}(\detd_V) @>{\bar\mu}>>{\Bbb V}(\detd_V)
\\ @VVV @VVV \\
\glv\times\grv @>{\mu}>> \grv
\endCD$$
\end{thm}

\begin{pf} Let us define $\glve(S)$ as the set of commutative
diagrams:
$$\CD {\Bbb V}(\detd_V) @>{\bar g}>>{\Bbb V}(\detd_V)
\\@VVV @VVV \\
\grv @>g>>\grv \endCD$$ for each $g\in\glv(S)$, and the
homomorphism
$\glve\to\glv$ given by $\bar g\mapsto g$. The rest of the proof
follows immediately from the fact that
$ H^0(\gr^0(V),\o_{\gr^0(V)})=k$ (remark {\ref{section-ogr}}) and
$g^*\detd_V\simeq\detd_V$ for every
$g\in\glv$.
\end{pf}

\begin{rem}\label{G-extension}
Let $G$ be a commutative subgroup of $\glv$ (a subfunctor of
commutative groups). The central extension of $\glv$ gives an
extension of $G$:
$$ 0\to {\Bbb G}_m\to\tilde G @>{\pi}>> G\to 0$$ and the commutator
of $\tilde G$:
$$\aligned
\tilde G\times\tilde G &\to\tilde G\\ (\tilde a,\tilde b)& \mapsto
\tilde a\tilde b\tilde a^{-1}\tilde b^{-1}
\endaligned$$ induces a pairing:
$$\aligned G\times G & @>{[\, , \,]}>> {\Bbb G}_m \\ (g_1,g_2)
&\mapsto [g_1,g_2]=
\tilde g_1\tilde g_2\tilde g_1^{-1}\tilde g_2^{-1}
\qquad \left(\tilde g_i\in\pi^{-1}(g_i)\right)
\endaligned$$ When $V$ is a local field or a ring of adeles, this
pairing will be of great importance in the study of arithmetic
problems because it is connected with the formulation of reprocity
laws.

The same construction of the extensions $\tilde G$ applied to the
Lie algebra of $G$ gives an extension of Lie algebras (taking the
points of $G$ with values in $k[x]/x^2$):
$$ 0\to {\Bbb G}_a=Lie({\Bbb G}_m)\to Lie(\tilde G) @>{d\pi}>>
Lie(G)\to 0$$ and a pairing:
$$\aligned Lie(G)\times Lie(G) & @>{R}>> {\Bbb G}_a \\ (D_1,D_2)
&\mapsto R((D_1,D_2))=[\tilde D_1,\tilde D_2]=
\tilde D_1\tilde D_2 - \tilde D_2\tilde D_1
\endaligned$$ ($\tilde D_i$ being a preimage of $D_i$).
\end{rem}

The pairing $R$ is an abstract generalization of the definition of
Tate [{\bf T}] of the residue pairing. There are several subgroups
of special relevance in the application of this theory to the study
of moduli problems and soliton equations. Firstly, we are concerned
with the algebraic analogue of the group $\Gamma$ ([{\bf SW}]
\S 2.3) of continuous maps $S^1\to{\Bbb C}^*$ acting as
multiplication operators over the Grassmannian. The main difference
between our definition of the group $\Gamma$ and the definitions
offered in the literature ([{\bf SW}], [{\bf PS}]) is that in the
algebro-geometric setting the elements
$\sum_{-\infty}^{+\infty}g_k\,z^k$ with infinite positive and
negative coefficients do not make sense as multiplication operators
over $k((z))$.

Let us now consider the case $V=k((t)), V^+=k[[t]]$. The main idea
for defining the algebraic analogue of the group
$\Gamma$ is to construct a ``scheme'' whose set of rational points
is precisely the multiplicative group
$k((z))^*$.

\begin{defn} The contravariant functor, $\kz$, over the category of
$k$-schemes with values in the category of commutative groups is
defined by:
$$S\rightsquigarrow \kz(S)= H^0(S,\o_S)((z))^*$$ Where for a
$k$-algebra $A$, $A((z))^*$ is the group of invertible elements of
the ring $A((z))=A[[z]][z^{-1}]$.
\end{defn}

\begin{lem}\label{v-loc-const} For each $k$-scheme $S$ and
$f\in\kz(S)$, the function:
$$\aligned S&\to \Z \\ s &\mapsto v_s(f)= \text{order of
$f_s\in k(s)((z))$}
\endaligned$$ is locally constant.
\end{lem}

\begin{pf} We can assume that $S=\sp(A)$, $A$ being a $k$-algebra.
Let
 $f=\sum_{i\ge n}^{}a_i\,z^i$ be an element of $A((z))^*$
($n\in\Z$). There then exists another element $g=\sum_{i\ge
-m}^{}b_i\,z^i$ ($m\in \Z$) such that $f\cdot g=1$. This
implies the following relations (from now on we assume $n=0$ to
simplify the calculations):
\beq
\begin{aligned}
 & 0=b_{-m}\, a_0 \\ & 0=b_{-m}\,a_1+b_{-m+1}\,a_0 \\ &
\dots \\ & 0=b_{-m}\,a_{m-1}+\dots+ b_{-1}\,a_0 \\ &
1=b_{-m}\,a_m+\dots+ b_0\,a_0
\end{aligned}\label{relations}
\end{equation}

Let us distinguish two cases:

\item{a)} $b_{-m}$ is not nilpotent in $A$: from {\ref{relations}}
we obtain:
$$b_{-m}\,a_0=b_{-m}^2\,a_1=\dots b_{-m}^m\,a_{m-1}=0$$ That is,
$a_0,\dots,a_{m-1}$ are equal zero in the ring
$A_{(b_{-m})}$ and for each $s\in\sp(A)-(b_{-m})_0$ one has
$b_{-m}(s)\,a_m(s)=1$ and therefore $v_{s}(f)=m$. We conclude by
proving that in this case $(b_{-m})_0$ is also an open subset of
$\sp(A)$:

\indent From the equations {\ref{relations}} we deduce:
$$\gathered  (b_{-m},a_{m-1},\dots,a_0)_0= (b_{-m})_0\cap
(a_{m-1})_0 \cap \dots\cap(a_0)_0=\emptyset\\
(b_{-m})_0\cup(a_i)_0=\sp(A)\quad,\quad i=0,\dots,n-1
\endgathered$$
\indent and hence:
$$(b_{-m})_0\cup\left(\cap_{i=0}^{n-1}(a_i)_0\right)=\sp(A)$$

\item{b)} Let us assume that $b_{-m},\dots,b_{-r-1}$ are nilpotent
elements of $A$ and that $b_{-r}$ is not nilpotent. The same
argument as in case a) proves that $v_s(f)$ is constant in the
closed subscheme $(b_{-r})_0$ and that its complementary in
$\sp(A)$ is $\cap_{i=0}^{r-1}(a_i)_0$, from which we conclude the
proof.
\end{pf}

\begin{cor} For an affine irreducible $k$-scheme $S=\sp(A)$ one has
that:
\begin{enumerate}
\item $v_s$ is a constant function over $S$,
\item
$$\left\{f\in A((z))^* \,\vert\, v(f)=n\right\}=
\left\{\gathered
\text{series }\,a_{n-r}\,z^r+\dots+a_n\,z^n+\dots\text{ such
that}\\ a_{n-r},\dots,a_{n-1}\text{ are nilpotent and }a_n\in A^*
\endgathered\right\}$$
\item If $A$ is also a reduced $k$-algebra:
$$A((z))^*=\coprod_{n\in \Z}\left\{\sum_{i\ge n} a_i\,z^i
\quad a_i\in A\text{ y }a_n\in A^*\right\}$$
\end{enumerate}
\end{cor}

\begin{pf} This is obvious from lemma {\ref{v-loc-const}}.
\end{pf}

\begin{thm} The subfunctor  $\kz_{red}$ of $\kz$ defined by:
$$S\rightsquigarrow \kz_{red}(S)=
\coprod_{n\in \Z}\left\{z^n+\sum_{i> n} a_i\,z^i \quad a_i\in
H^0(S,\o_S)\right\}$$ is representable by a group $k$-scheme whose
connected component of the origin will be denote by $\Gamma_+$.
\end{thm}

\begin{pf} It suffices to observe that the functor:
$$S\rightsquigarrow \left\{z^n+\sum_{i> n} a_i\,z^i \quad a_i\in
H^0(S,\o_S)\right\}$$ is representable by the scheme:
$$\sp(\limil{l} k[x_1,\dots,x_l])=\limpl{l}\A^l_k$$ and the
group law is given  by the multiplication of series.
\end{pf}

\begin{thm} Let $\kz_{nil}$ be the subfunctor of $\kz$ defined by:
$$S\rightsquigarrow \kz_{nil}(S)=
\coprod_{n>0}\left\{\gathered
\text{ finite series }\,a_n\,z^{-n}+\dots+a_1\,z^{-1}+1\text{ such
that } \\  a_i\in H^0(S,\o_S)\text{ are nilpotent and $n$ arbitrary}
\endgathered\right\}$$ There exists a formal $k$-scheme $\Gamma_-$
representing
$\kz_{nil}$, that is:
$$\hom_{\text{for-sch}}(S,\Gamma_-)=\kz_{nil}(S)$$ for every
$k$-scheme $S$.
\end{thm}

\begin{pf} Let us define the ring of ``infinite'' formal series in
infinite variables (which is different from the ring of formal
series in infinite variables) by:
$$k\{\{x_1,\dots\}\}=\underset n\limp k[[x_1,\dots,x_n]]$$ the
morphisms of the projective system being:
$$\aligned k[[x_1,\dots,x_{n+1}]]&\to k[[x_1,\dots,x_n]] \\
x_i&\mapsto x_i\qquad\text{for }i=1,\dots,n-1 \\ x_{n+1}&\mapsto 0
\endaligned$$ Note that:
$$k\{\{x_1,\dots\}\}=\underset
n\limp{k[x_1,\dots,x_n]}/{(x_1,\dots,x_n)^n}$$ It is therefore an
admissible linearly topological ring  ([{\bf EGA}]~{\bf 0}.7.1) and
there therefore exists its formal spectrum
$\sf(k\{\{x_1,\dots\}\})$. Let us denote by $J_n$ the kernel of the
natural projection
$k\{\{x_1,\dots\}\}\to {k[x_1,\dots,x_n]}/{(x_1,\dots,x_n)^n}$ and
 $J=\limp (x_1,\dots,x_n)$.

Let us now prove that $\Gamma_-=\sf(k\{\{x_1,\dots\}\})$:

\noindent For every $k$-scheme $S$, considering over
$H^0(S,\o_S)$ the discrete topology, we have:
$$\aligned\hom_{\text{for-sch}}(S,\Gamma_-)=&
\hom_{\text{cont-$k$-alg}}((k\{\{x_1,\dots\}\},H^0(S,\o_S))=
\\=&\left\{\gathered
f\in\hom_{\text{$k$-alg}}(k\{\{x_1,\dots\}\},H^0(S,\o_S))
\text{ such that} \\
\text{there exists $n\in{\Bbb N}$ satisfying }J_n\subseteq
f^{-1}((0))\endgathered\right\}
\endaligned$$  However the condition $J_n\subseteq f^{-1}((0))$ is
equivalent to saying that $f(x_1),\dots,f(x_n)$ are nilpotent and
$f(x_i)=0$ for
$i>n$, from which one concludes the proof.
\end{pf}

\begin{rem} Note that $\Gamma_-$ is the inductive limit in the
category of formal schemes ([{\bf EGA}]~{\bf I}.10.6.3) of the
schemes which represent the subfunctors:
$$S\rightsquigarrow \Gamma^n_-(S)=\left\{\gathered
\,a_n\,z^{-n}+\dots+a_1\,z^{-1}+1\text{ such that}\\ a_i\in
H^0(S,\o_S)\text{ and the ${\text{n}}^{\text{th}}$ power} \\
\text{ of the ideal $(a_1,\dots,a_n)$ is zero}
\endgathered\right\}$$
\end{rem}

\begin{rem}{\bf Group laws of $\Gamma_+$ and $\Gamma_-$} The group
law of $\Gamma_+=\sp(k[x_1,\dots])$ is given by:
$$\aligned  k[x_1,\dots]&\to k[x_1,\dots]\otimes_k k[x_1,\dots] \\
x_i &\mapsto x_i\otimes 1+\sum_{j+k=i}x_j\otimes x_k+1\otimes x_i
\endaligned$$ The group law of $\Gamma_-=\sf k\{\{x_1,\dots\}\}$ is
given by:
$$\aligned k\{\{x_1,\dots\}\}&\to k\{\{x_1,\dots\}\}\hat\otimes_k
k\{\{x_1,\dots\}\} \\ x_i &\mapsto x_i\otimes
1+\sum_{j+k=i}x_j\otimes x_k+1\otimes x_i
\endaligned$$
\end{rem}

Let be $\kz_0$ be the connected component of the origin in the
functor of groups $\kz$.

\begin{thm} The natural morphism of functors of groups over the
category of
$k$-schemes:
$$\fu{\Gamma_-}\times\fu{{\Bbb G}_m}\times\fu{\Gamma_+}\to\kz$$ is
injective and for $char(k)=0$ gives an isomorphism with
$\kz_0$. $\kz_0$ is therefore representable by the (formal)
$k$-scheme:
$$\Gamma=\Gamma_-\times{\Bbb G}_m\times\Gamma_+$$
\end{thm}

\begin{pf} The morphism from $\fu{{\Bbb G}_m}$ to $\kz$ is the one
induced by the natural inclusion $H^0(S,\o_S)^*\hookrightarrow
H^0(S,\o_S)((z))^*$.

The injectivity of $\fu\Gamma\hookrightarrow\kz$ follows from the
fact that $\Gamma_-\cap\Gamma_+=\{1\}$. The rest of the proof is
trivial from the above results and from the properties of the
exponential map we shall see below.
\end{pf}

\begin{rem} Our group scheme $\Gamma$ is the algebraic analogue of
the group $\Gamma$ of Segal-Wilson [{\bf SW}]. Note that the
indexes ``-'' and ``+'' do not coincide with the Segal-Wilson
notations. Replacing $k((z))$, by $k((z^{-1}))$ we obtain the same
notation as in the paper of Segal-Wilson.
\end{rem}

Let us define the exponential maps for the groups $\Gamma_-$ and
$\Gamma_+$. Let $\A_n$ be the $n$ dimensional affine space
over $\spk$ with the additive group law, and $\hat\A_n$ the
formal group obtained as the completion of $\A_n$ at the
origin. We define $\hat\A_\infty$ as the formal group
$\limil{n}\hat\A_n$. Obviously $\hat\A_\infty$ is the
formal scheme:
$$\hat\A_\infty=\sf k\{\{y_1,\dots\}\}$$ with group law:
$$\aligned k\{\{y_1,\dots\}\}&\to  k\{\{y_1,\dots\}\}\hat\otimes_k
k\{\{y_1,\dots\}\} \\ y_i& \longmapsto y_i\otimes 1+1\otimes y_i
\endaligned$$

\begin{defn} If the characteristic of $k$ is zero, the exponential
map for
$\Gamma_-$ is the following isomorphism of formal group schemes:
$$\aligned
\hat\A_\infty &@>{\exp}>> \Gamma_- \\
\{a_i\}_{i>0} &\mapsto \exp(\sum_{i>0}a_i\,z^{-i})
\endaligned$$ This is the morphism induced by the ring homomorphism:
$$\aligned
 k\{\{x_1,\dots\}\}& @>{\qquad\exp^*\qquad}>> k\{\{y_1,\dots\}\}\\
x_i &\mapsto \text{ coefficient of $z^{-i}$ in the series }
\exp(\sum_{j>0}y_j\,z^{-j})\endaligned$$
\end{defn}

\begin{defn}\label{exp-gamma-minus} If the characteristic of $k$ is
$p>0$, the exponential map for
$\Gamma_-$ is the following isomorphism of formal schemes:
$$\aligned \hat\A_\infty &\to \Gamma_- \\
\{a_i\}_{i>0}&\mapsto \prod_{i>0}(1-a_i\,z^{-i})\endaligned$$ which
is the morphism induced by the ring homomorphism:
$$\aligned k\{\{x_1,\dots\}\}&@>\exp^*>> k\{\{y_1,\dots\}\}\\ x_i
&\mapsto \text{ coefficient of $z^{-i}$ in the series }
\prod_{i>0}(1-a_i\,z^{-i})\endaligned$$
\end{defn}

Note that this latter exponential map is not a isomorphism of
groups. Considering over $\hat \A_\infty$ the law group
induced by the isomorphism, $\exp$, of formal schemes, we obtain
the Witt formal group law.

Analogously, we define the exponential maps for the group
$\Gamma_+$:

\begin{defn} Let $\A^\infty$ be the group scheme over $k$
defined by
$\limpl{n} \A_n$ (where $\A_{n+1}=\sp
k[x_1,\dots,x_{n+1}]\to \A_n=\sp k[x_1,\dots,x_n]$ is the
morphism defined by forgetting the last coordinate) with its
additive group law. The exponential map when $char(k)=0$ is the
isomorphism of group schemes:
$$\aligned  \A^\infty &\to \Gamma_+ \\
\{a_i\}_{i>0}&\mapsto \exp(\sum_{i>0}a_i\,z^i)\endaligned$$ If
$char(k)=p\ne 0$, the exponential map is the isomorphism of schemes:
$$\aligned \A^\infty &\to \Gamma_+ \\
\{a_i\}_{i>0}&\mapsto \prod_{i>0}(1-a_i\,z^i)\endaligned$$ which is
not a morphism of groups.

(See {\rm [{\bf B}]} for the connection of these definitions and the
Cartier-Dieudonn\'e theory).
\end{defn}

It should be noted that the formal group scheme $\Gamma_-$ has
properties formally analogous to the Jacobians of the algebraic
curves: one can define formal Abel maps and prove formal analogues
of the Albanese property of the Jacobians of smooth curves (see
{\rm [{\bf KSU},{\bf C}]}).

Let $\hat C=\sf(k[[t]])$ be a formal curve. We define the Abel
morphism of degree $1$ as the morphism of formal schemes:
$$\phi_1: \hat C\to \Gamma_-$$ given by
$\phi_1(t)=(1-\frac{t}z)^{-1}=1+\sum_{i>0}^{}\frac{t^i}{z^i}$; that
is, the morphism induced by the ring homomorphism:
$$\aligned k\{\{x_1,\dots\}\}&\to k[[t]]\\ x_i\,&\mapsto
t^i\endaligned$$

Note that the Abel morphism is the algebro-geometric version of the
function $q_\xi(z)$ used by Segal and Wilson ([{\bf SW}] page~32)
to study the Baker function.

Let us explain further why we call $\phi_1$ the ``Abel morphism''
of degree 1. If $char(k)=0$, composing
$\phi_1$ with the inverse of the exponential map, we have:
$$\bar\phi_1:\hat C@>{\phi_1}>> \Gamma_-@>>{\exp^{-1}}>>
\hat\A_\infty$$ and since
$(1-\frac{t}z)^{-1}=\exp(\sum_{i>0}^{}\frac{t^i}{i\,z^i})$ (see
[{\bf SW}] page~33), $\bar\phi_1$ is the morphism defined by the
ring homomorphism:
$$\aligned k\{\{y_1,\dots\}\} &\to k[[t]]\\ y_i&\mapsto \frac{t^i}i
\endaligned$$ or in terms of the functor of points:
$$\aligned\hat C&@>{\bar\phi_1}>>\hat\A_\infty\\ t&\mapsto
\{t,\frac{t^2}2,\frac{t^3}3,\dots\}
\endaligned$$ Observe that given the basis $\omega_i=t^i\,dt$ of the
differentials $\Omega_{\hat C}=k[[t]]dt$, $\bar\phi_1$ can be
interpreted as the morphism defined by the ``abelian integrals''
over the formal curve:
$$\bar\phi_1(t)=\left(
\int\omega_0,\int\omega_1,\dots,\int\omega_i,\dots
\right)$$ which coincides precisely with the local equations of the
Abel morphism for smooth algebraic curves over the field of complex
numbers. In general, for each integer number $n>0$, we define the
Abel morphism of degree $n$ as the morphism of formal schemes:
$$\bar\phi_n:\hat C\times\overset n\dots\times\hat C=\hat
C^n\to\Gamma_-$$ given by
$\bar\phi_n(t_1,\dots,t_n)=
\prod_{i=1}^n\left(1-\frac{t_i}z\right)^{-1}$; that is, the
morphism induced by the ring homomorphism:
$$\aligned k\{\{x_1,\dots\}\}&\to  k[[t_1]]\hat\otimes\overset
n\dots\hat\otimes k[[t_n]]\\ x_i\,&\mapsto \text{ coefficient of
$z^{-i}$ in the series }
\prod_{i=1}^n(1-\frac{t_i}{z})^{-1}\endaligned$$ Note that
$\bar\phi_n$ factorizes through a morphism,
$\phi_n$ from the $n^{\text{th}}$-symmetric product of $\hat C$ to
$\Gamma_-$, which is the true Abel morphism; moreover $\phi_n$ is
an immersion.

\begin{thm} $(\Gamma_-,\phi_1)$ satisfies the Albanese property for
$\hat C$; that is, every morphism $\psi:\hat C\to X$ in a
commutative group scheme (which sends the unique rational point of
$\hat C$ to the $0\in X$) factors through the Abel morphism and a
homomorphism of groups $\Gamma_-\to X$.
\end{thm}

\begin{pf} Let $\psi:\hat C\to X$ be a morphism from the formal
scheme
$\hat C$ to a group scheme $X$ such that $\psi(\text{rational
point})=0$. For each $n>0$, one constructs a morphism:
$$\hat C^n @>{\bar\psi_n}>> X$$ which is the composition of
$\psi\times\dots\times\psi: \hat C\times\overset n\dots\times\hat C
\to X\times\overset n\dots\times X$ and the addition morphism
$X\times\overset n\dots\times X\to X$. Observe that
$\bar \psi_n$ factors through a morphism:
$$S^n\hat C@>{\psi_n}>> X$$ and bearing in mind that
$\Gamma_-=\limil{n} S^n\hat C$ (as formal group schemes) we
conclude the proof of the existence of a homomorphism of groups
$\bar\psi:\Gamma_-\to X$ satisfying the desired condition.
\end{pf}

%%%%%%%%%%%%%%%%%%%%%%%%%%%%%%%%%%%%%%%%%%%%%%%%%%%%%%%%%%%%%
%%%
\section{$\tau$-functions and Baker functions}

This section is devoted to algebraically defining the
$\tau$-functions and the Baker functions over an arbitrary base
field $k$.

Following on with the analogy between the groups $\Gamma$ and
$\Gamma_-$ and the Jacobian of the smooth algebraic curves, we
shall make the well known constructions for the jacobians of the
algebraic curves for the formal curve $\hat C$ and the group
$\Gamma$: Poincar\'e bundle over the dual jacobian and the
universal line bundle over the jacobian. In the formal case these
constructions are essentially equivalent to defining the
$\tau$-functions and the Baker functions.

Using the notations of section {\ref{aut-grass}}, let us consider
the Grassmannian $\grv$ of $V=k((z))$  and the group
$$\Gamma=\Gamma_-\times {\Bbb G}_m \times \Gamma_+$$  acting on
$\grv$ by homotheties.

As we have shown in {\ref{G-extension}}, there exists a central
extension of
$\Gamma$:

$$0\to {\Bbb G}_m \to \widetilde \Gamma \to \Gamma \to 0$$
given by a pairing:
\beq\quad \Gamma\times \Gamma \to {\Bbb G}_m
\label{pairing}\end{equation}

\begin{prop}
The extension $\tilde\Gamma_+$ of $\Gamma_+$ is trivial.
\end{prop}

\begin{pf}
We will construct a section $s$ (as groups) of
$\tilde\Gamma_+\to\Gamma_+$; that is, for an element
$g\in\Gamma_+$ we give $s(g)\in\tilde\Gamma_+$ such that
$s$ is a morphism of groups.

Denote by $\mu:\Gamma\times\grv\to\grv$ the action of
$\Gamma$ on $\grv$ and by $\mu_g$ the automorphism of
$\grv$ induced by the homothety $\cdot g:V\to V$ for
$g\in\Gamma$.

Fix $g\in\Gamma_+$. Observe that there exists a
quasi-isomorphism of complexes:
$$\CD \L @>>> V/V^+ \\ @V{\cdot g}VV>  @VV{\cdot g}V \\
\mu_g^*(\L) @>>> V/V^+ \endCD$$
since $g\cdot V^+\simeq V^+$. We have thus an isomorphism
$\detd_V\simeq \mu_g^*\detd_V$ in a canonical way, and
hence a well-defined element $s(g)\in\tilde\Gamma_+$.
Since this construction is canonical and
$\mu_{g'}\circ\mu_g=\mu_{g'\cdot g}$ it follows easily that
$s(g')\cdot s(g)=s(g'\cdot g)$.
\end{pf}

\begin{prop}
For a rational point $U\in\grv$, let $\mu_U$ be the
morphism $\Gamma\times\{U\}\to\grv$ induced by $\mu$.
Then, the line bundle $\mu_U^*\detd_V\vert\Gamma_-$ is
trivial, and the extension $\tilde\Gamma_-$ is thus
trivial.
\end{prop}

\begin{pf}
Assume $U\in F_{V^+}$ (the general case is anologous).
It is no difficult to obtain the following equality for
$g\in\Gamma_-$:
$$(\mu_U^*\Omega_+)(g)\,=\,\Omega_+(g\cdot U)\,=\,
\Omega_+(U)+\sum_S \chi_S(g)\cdot\Omega_S(U)$$
where the sum is taken over the set of Young diagrams and
$\chi_S$ is the Schur polynomial (in the coefficients of
$g$) corresponding to $S$. Since $\Omega_+(U)\neq0$ and
the coefficients of $g$ are nilpotents, it follows that
$\mu_U^*\Omega_+$ is a no-where vanishing section of
$\mu_U^*\detd_V$, and this bundle is therefore trivial.

Observe now that since $\tilde\Gamma_-$ can be thought as
the sheaf of automorphisms of $\mu_U^*\detd_V$ one has
that $\tilde\Gamma_-$ is a trivial extension.
\end{pf}

\begin{cor}
The restrictions of the pairing {\ref{pairing}} to the subgroups
$\Gamma_-$  and $\Gamma_+$ are trivial.
\end{cor}

We define the Poincar\'e bundle over $\Gamma\times \grv$ as the
invertible sheaf:
$${\frak P}=\mu^*\det_V^*$$

For each point $U\in \grv$, let us define the Poincar\'e bundle over
$\Gamma\times \Gamma$ associated with $U$ by:
$${\frak P}_U=(1\times \mu_U)^*{\frak P}=m^*({\mu_U}^*\det_V^*)$$
where $m:\Gamma\times\Gamma\to\Gamma$ is the group law.

The sheaf of $\tau$-functions of a point $U\in \grv$, $\widetilde
{\L_\tau}(U)$, is the invertible sheaf over $\Gamma\times \{U\}$
defined by:
$$\widetilde{\L_\tau}(U)={\frak P}\vert_{\Gamma\times\{U\}}$$

Let us note that the sheaf $\widetilde{\L_\tau}(U)$ is defined for
arbitrary points of the Grassmannian and not only for geometric
points.

The restriction homomorphism induces the following homomorphism
between global sections:
\beq H^0\left(\Gamma\times \grv, \mu^*\det_V^*\right)\to
 H^0\left(\Gamma\times \{U\},\widetilde {\L_\tau}(U)  \right)
\label{restriction}\end{equation}

\begin{defn}
The $\tau$-function of the point $U$ over $\Gamma$ is
defined as the image $\widetilde {\tau }_U$ of the section
$\mu^*\Omega_+$ by the homomorphism {\ref{restriction}} ($\Omega_+$
being the global section defined in {\ref{global-section}}).
\end{defn}

Obviously $\widetilde {\tau }_U$ is not a function over
$\Gamma\times
\{U\}$ since the invertible sheaf $\widetilde {\L_\tau}(U)$ is not
trivial.

The algebraic analogue of the $\tau$-function defined by M. and Y.
Sato, Segal and Wilson ([{\bf SS}], [{\bf SW}]) is obtained by
restricting the invertible sheaf $\widetilde {\L_\tau}(U)$ to the
formal subgroup $\Gamma_-\subset\Gamma$.

To see this, fix a rational point $U\in\grv$ and define:
$$\L_\tau(U)=\widetilde{\L_\tau}(U)\vert_{\Gamma_-\times\{U\}}$$
which is a trivial invertible sheaf over $\Gamma_-$. To
obtain a trivialization of ${\L_\tau}(U)$ which will allow us to
identify global sections with functions over $\Gamma_-$ we must fix
a global section of ${\L_\tau}(U)$ without zeroes in $\Gamma_-$.

Recall that $\tilde\Gamma_-$ is a trivial extension of $\Gamma_-$
and it has therefore a section $s$. It follows that the group
$\Gamma_-$ acts on $\L_\tau(U)$ (through $s$) and on $\Gamma_-$ by
translations. One has easily that the morphism:
$${\mathbb V}(\L_\tau(U)^*)\to\Gamma_-$$
 is equivariant with respect to these actions.

Note now that:
$$\hom_{\Gamma_-\text{-equiv}}
\left(\Gamma_-,{\mathbb V}(\L_\tau(U)^*)\right)\subseteq
\hom_{\Gamma_-\text{-esq}}
\left(\Gamma_-,{\mathbb V}(\L_\tau(U)^*)\right)=
H^0(\Gamma_-,\L_\tau(U))$$
Let $\delta$ be an non-zero element in the fibre of
${\mathbb V}(\L_\tau(U)^*)$ over the point $1$ of $\Gamma_-$ (1
being the identity of $\Gamma_-$). Let  $\sigma_0$ be the unique
morphism $\Gamma_-\to{\mathbb V}(\L_\tau(U)^*)$
$\Gamma_-$-equivariant such that $\sigma_0(1)=\delta$, and denote
again by $\sigma_0$ the corresponding section of $\L_\tau(U)$.

Observe that $\sigma_0$ is a constant section and since it has no
zeros it gives a trivialization of $\L_\tau(U)$. Through this
trivilization, the global section of ${\L_\tau}(U)$ defined by
$\tilde\tau_U$ is identified with the function
$\tau_U\in\o(\Gamma_-)= k\{\{x_1,\dots\}\}$ given by
Segal-Wilson [{\bf SW}]:
$${\tau}_U(g)=\frac{\tilde\tau_U(g)}{\sigma_0(g)}=
\frac{\mu^*\Omega_+(g)}{\sigma_0(g)}=\frac{\Omega_+(gU)}{\delta}$$
Finally, if $U\in F_{V^+}$ then one can choose $\delta=\Omega_+(U)$.

Observe that the $\tau$-function ${\tau}_U$ is not a series of
infinite variables but an element of the ring $k\{\{x_1,
\dots\}\}$.

The subgroup $\Gamma_+$ of $\Gamma$ acts freely over $\grv$.
Accordingly the orbits of the rational points of $\grv$ under the
action of $\Gamma_+$ are isomorphic, as schemes, to $\Gamma_+$.

Let $X$ be the orbit of $V^-=z^{-1}\cdot k[z^{-1}]\subset V$ under
$\Gamma_+$. The restrictions of $\det_V$ and $\detd_V$ to $X$ are
trivial invertible sheaves. Bearing in mind that the points of
$X$ are $k$-vector subspaces of $V$ whose intersection with $V^+$ is
zero, one has that the section $\Omega_+$ of $\det_V^* $ defines a
canonical trivialization of $\detd_V$ over $X$.

\begin{thm}
The restriction homomorphism $\detd_V
\to\detd_V\vert_X$ induces a homomorphism between global sections:
$$ B: H^0\left(\grv,\detd_V \right)\longrightarrow
 H^0\left(X,\detd_V\vert_X \right)\simeq
\o(\Gamma_+)=k[x_1,\dots]$$
 which is an isomorphism between the $k$-vector subspace $\Omega(S)$
defined in {\ref{global-section}} and $\o(\Gamma_+)$. The
isomorphism
$ H^0(X,\detd_V\vert_X)\overset \sim\rightarrow \o(\Gamma_+)$ is
the isomorphism induced by the trivialization defined by
$\Omega_+$.

In the literature, the isomorphism $B:\Omega(S)\overset \sim
\longrightarrow
\o(\Gamma_+)$ is usually called the bosonization isomorphism.
\end{thm}

\begin{pf}
All one has to prove is that $B(\Omega_S)=F_S(x)$,
$\Omega_S$ being the Pl\"ucker sections of $\detd_V$ defined in
{\ref{global-section}} and $F_S(x_1,x_2,\dots)$ being the Schur
functions. Proof of the identity $B(\Omega_S)=F_S(x)$ is
essentially the same as in the complex analytic case; see [{\bf
SW}] and [{\bf PS}].

In some of the literature, the $\tau$ function of a point
$U\in\grv$ is defined as the Pl\"ucker coordinates of the point
$U$. Let us therefore explain in which sense both definitions are
equivalent.

The canonical homomorphism:
$$H^0(\detd_V)\otimes\o_{\grv}\longrightarrow\detd_V\to 0$$
induces a homomorphism:
$$\det_V=\det_V^{**}\overset{\bar\tau}\hookrightarrow
H^0(\detd_V)^*\otimes\o_{\grv}$$
\end{pf}

\begin{defn}
Given a point $\widetilde U\in\det_V$ in the fibre of
$U\in\grv$, the $\bar\tau$-function of $\widetilde U$ is defined as
the element $\bar \tau(U)\in H^0\left (\detd_V\right)^*\otimes k(U)$
($k(U)$ being the residual field of $U$). This is essentially the
definition of $\tau$-functions given in the papers of M. and Y.
Sato, Arbarello and De Concini, and Kawamoto and others ({\rm [{\bf
SS}], [{\bf AD}], [{\bf KNTY}]}).
\end{defn}

\begin{lem}
There exists an isomorphism of $k$-vector spaces:
$$\o(\Gamma_+)^*\to \o(\Gamma_-)$$
\end{lem}

\begin{pf}
Recall that $\o(\Gamma_+)=k[x_1,\dots]$ and that
$\o(\Gamma_-)=k\{\{x_1,x_2,\dots\}\}$. Now think that $x_i$ is the
$i$-symmetric function of other variables, say $t_1,t_2,\dots$. It
is known that the Schur polynomials $\{F_S\}$ (where $S$ is a
partition) of the $t$'s are polynomials in the $x$'s and are in
fact a basis of the $k$-vector space $k[x_1,\dots]$. The isomorphism
is the induced by the pairing:
$$\aligned \o(\Gamma_+)\times \o(\Gamma_-) &\longrightarrow k \\
(F_S,F_{S'}) & \longmapsto \delta_{S,S'}\endaligned$$
(see [{\bf Mc}]).
\end{pf}

The composition of the homomorphism $B^*$ (the dual homomorphism of
$B$) and the isomorphism of the above lemma gives an homomorphism:
$$\tilde B^*:\o(\Gamma_+)^*=k\{\{x_1,x_2,\dots\}\}
\longrightarrow H^0\left (\detd_V\right)^*$$

The connection between $\tau_U$ and $\bar \tau (\widetilde U)$ is
the following:
$$\tilde B^*(\tau_U)=\lambda \cdot(\bar \tau (\widetilde U))$$
$\lambda$ being a non-zero constant. (Of course, if $U$ is
not rational but a point with values in a scheme $S$, $\lambda
\in H^0\left(S,\o_S\right)^*$).

The connection of the $\tau$-functions with autoduality (in the
sense of group schemes) properties of the group $\Gamma=
\Gamma_-\times {\Bbb G}_m\times \Gamma_+$ implicit in the above
discussion, is studied with detail in [{\bf C},{\bf P}].
L. Breen in [{\bf B2}] outlines also some of these properties from
another point of view.

Once we have algebraically defined the $\tau$-functions, we can
define the Baker functions using formula~5.14. of [{\bf SW}]; this
is the procedure used by several authors. However, we prefer to
continue with the analogy with the classical theory of curves and
jacobians and define the Baker functions as a formal analogue of
the universal invertible sheaf of the Jacobian.

Let us consider the composition of morphisms:
$$\tilde\beta:\hat C\times\Gamma\times\grv
\overset {\phi\times Id}\longrightarrow
\Gamma\times\Gamma\times\grv\overset{m\times Id}
\longrightarrow\Gamma\times\grv$$

$\phi:\hat C=\sf k[[z]] \to \Gamma$ being the Abel morphism (taking
values in $\Gamma_-\subset \Gamma$) and
$m:\Gamma\times \Gamma\to\Gamma$ the group law.

\begin{defn}
The sheaf of Baker-Akhiezer functions is the
invertible sheaf over $\hat C\times \Gamma\times \grv$ defined by:
$$\widetilde {\L_B}=(\phi\times Id)^*(m\times Id)^*{\frak P}$$

Let us define the sheaf of Baker functions at a point $U\in \grv$ as
the invertible sheaf:
$$\widetilde {\L_B}(U)=\widetilde {\L_B}\vert_{\hat C\times
\Gamma\times \{U\}}=\tilde {\beta_U}^*\widetilde {\L_\tau}(U)$$
(where $\tilde {\beta_U}^*$ is the following homomorphism between
global sections:
$$ H^0(\Gamma\times \{U\},\widetilde {\L_\tau}(U))
\overset{\widetilde{\beta_U}^*}\longrightarrow
 H^0(\hat C\times\Gamma\times \{U\},\widetilde {\L_B}(U))$$
\end{defn}

By the definitions, $\widetilde {\L_B}(U)\vert_{\hat C\times
\Gamma_-\times
\{U\}}={\L_B}(U)$  is a trivial invertible sheaf over $\hat C \times
\Gamma_-$.

Observe that for each element $u\in\Gamma_-(S)\subseteq
\fu{k((z))^*}(S)=H^0(S,\o_S)((z))^*$ we can define a fractionary
ideal of the formal curve $\w {C_S}$ by:
$$I_u=u\cdot\o_S((z))$$ in such a way that we can interpret the
formal group $\Gamma_-$ as a kind of Picard scheme over the formal
curve. The universal element of $\Gamma_-$ is the invertible
element of $\kz(\Gamma_-)$ given by:
$$v=1+\underset {i \geq 1}\sum x_i\,z^{-i}\in k((z))\hat \otimes
k\{\{x_{ 1},x_{ 2},\dots\}\}$$
This universal element will be the formal analogue of the universal
invertible sheaf for the formal curve $\hat C$.

\begin{defn}
The Baker function of a point $U\in \grv$ is
 $\psi_U=v^{-1}\cdot \beta^*_U(\tilde\tau_U)$, where
$$\beta_U^*: H^0\left( \Gamma\times \{U\},\widetilde
{\L_\tau}(U)\right)\longrightarrow
 H^0\left(\hat C\times
\Gamma\times \{U\},\widetilde {\L_B}(U)\right)$$
is the homomorphism induced by $\widetilde{\beta_U}^*$.
\end{defn}

Observe that the Baker function of $V^-=z^{-1}\,k[z^{-1}]$ is the
universal invertible element $v^{-1}$.

Note that, analogously to the case of $\tau$-function, we can choose
a trivialization of $\widetilde {\L_B}(U)$ over $\hat
C\times\Gamma_-\times \{U\}$ in such a way that the
function asociated to the section $v^{-1}\cdot \beta^*_U(\tau_U)$ is:
\beq\psi_U(z,g)=v^{-1}\cdot \frac{\tau_U\left(g\cdot
\phi_1\right)}{\tau_U(g)}
\label{baker-expr}\end{equation}
which is the classical expression for  the Baker function.

When the characteristic of the base field $k$ is zero, we can
identify $\Gamma_-$ with the additive group scheme
 $\hat{\A}_{\infty}$ through the exponential and expression
{\ref{baker-expr}} is the classical expression for the Baker
functions ([{\bf SW}]~5.16):
$$\psi_U(z , t)=\left(  \frac{\tau_U(t+[z])}{\tau_U(t)
}\right)\cdot \exp(-\sum t_i\,z^{-i})$$  where
$[z]=(z,\frac{1}{2}z^2,\frac{1}{3} z^3,\dots)$ and
$t=(t_1,t_2,\dots)$ and $v=\exp(\sum t_i\,z^{-i})$ through the
exponential map.

For the general case, we obtain explicit expressions for $\psi_U$
as a function over
$\hat C\times \hat\A_{\infty}$ but considering in
$\hat\A_{\infty}$ the group law induced by the exponential
{\ref{exp-gamma-minus}}:
$$\psi_U(z,g)=v(z,g)^{-1}\cdot \frac{\tau_U\left(t*
\phi(z)\right)}{\tau_U(t) }$$ ($*$ being the group law of
$\hat\A_\infty$).

The classical properties characterizing the Baker functions (for
example proposition~5.1 of [{\bf SW}]) can be immediately
generalized for the Baker functions over arbitrary fields.

\begin{rem}
Note that our definitions of
$\tau$-functions and Baker functions are valid over arbitrary base
fields and that can be generalized for $\Z$. One then has the notion
of $\tau$-function and Baker functions for families of elements of
$\grv$ and, if we consider the Grassmannian of $\Z((z))$ one then
has $\tau$-functions and Baker functions of the rational points of
$\gr\left(\Z((z))\right)$  and the geometric properties studied in
this paper have a translation into arithmetic properties of the
elements of
$\gr\left(\Z((z))\right)$. The results stated by Anderson in
[{\bf A}] are a particular case of a much more general setting
valid not only for
$p$-adic fields but also for arbitrary global field numbers. Our
future aims are to study the arithmetic properties related to these
constructions.
\end{rem}

%%%%%%%%%%%%%%%%%%%%%%%%%%%%%%%%%%%%%%%%%%%%%%%%%%%%%%%%%%

%%%%%%%%%%%%%%%%%%%%%%%%%%%%%%%%%%%%%%%%%%%%%%

\vskip1.5truecm { \'Alvarez V\'azquez, Arturo}\newline {\it
e-mail: } aalvarez@@gugu.usal.es
\vskip.3truecm { Mu\~noz Porras, Jos\'e Mar\'{\i}a}\newline
{\it e-mail: } jmp@@gugu.usal.es
\vskip.3truecm { Plaza Mart\'{\i}n, Francisco Jos\'e}\newline
{\it e-mail: } fplaza@@gugu.usal.es\hfill{\tiny Use this one to
contact.}

\end{document}